\documentclass[12pt]{article}
\usepackage{amsmath,amssymb,amsfonts,amsthm}
\usepackage{yhmath,relsize}
\usepackage{eucal,oldgerm}
\usepackage[dvips]{graphics}
\usepackage{pstricks,pst-plot,pst-node,pst-coil}
\newpsobject{showgrid}{psgrid}{subgriddiv=1,griddots=5,gridlabels=0pt}
\usepackage[arrow,matrix,tips,2cell]{xy}

\newcommand{\rang}{\right\rangle}

\newcommand{\com}{{\mathbb C}}

\newcommand{\proj}{\mathbf{P}}

\newcommand{\PP}{{\mathsf{P}}}

\newcommand\FF{\mathbb F}
\newcommand\ZZ{\mathsf Z}
\newcommand{\oh}{{\mathcal O}}
\newcommand{\T}{{\mathbf{T}}}
\newcommand{\OO}{{\mathcal{O}}}

\newcommand\II{{\mathcal I}}
\newcommand\bull{{\scriptscriptstyle\bullet}}
\newcommand\udot{^\bull}

\newcommand{\LL}{{\mathbb{L}}}

\newcommand{\Q}{\mathbb{Q}}
\newcommand{\Z}{\mathbb{Z}}

\newcommand{\rarr}{\rightarrow}

\newcommand{\bZ}{\mathsf{Z}}
\newcommand{\oM}{\overline{M}}

\DeclareMathOperator{\Hilb}{Hilb}

\DeclareMathOperator{\Ext}{Ext}

\newtheorem{Theorem}{Theorem}

\newtheorem{definition}{Definition}

\newtheorem{Conjecture}{Conjecture}

\begin{document}
\title{Descendents for stable pairs on 3-folds}
\author{R. Pandharipande\vspace{12pt}\\\vspace{5pt}
{\em {\footnotesize{Dedicated to Simon Donaldson on the occasion of his $60^{th}$ birthday}}}}

\date{March 2017}
\maketitle

\begin{abstract}
We survey here the construction and the basic properties
of descendent invariants in the theory of
stable pairs on nonsingular projective 3-folds. The main topics
covered are the rationality of the 
generating series, the functional equation,
the Gromov-Witten/Pairs 
correspondence for descendents, the Virasoro constraints,
and the connection to the virtual fundamental class of the stable
pairs moduli space
 in algebraic cobordism.
In all of these directions, the proven results constitute 
only a small part of the conjectural framework. A central
goal of the article is 
to introduce the open questions as simply and directly as possible.
\end{abstract}

\maketitle

\setcounter{tocdepth}{1} 
\tableofcontents


\setcounter{section}{-1}
\section{Introduction}

\subsection{Moduli space of stable pairs}
Let $X$ be a nonsingular projective $3$-fold.
The moduli of curves in $X$ can be approached in several different
ways.{\footnote{For a discussion of the different approaches, see \cite{rp13}.}} For an algebraic 
geometer, perhaps the most straightforward is
the Hilbert scheme of subcurves of $X$. 
The moduli of stable pairs is  
closely related to the Hilbert scheme, 
but is geometrically much more efficient. 
While the definition of a  stable pair takes some time to understand, 
the advantages of the moduli theory 
more than justify the effort.

\begin{definition} \label{sp}
A {\em{stable pair}} $(F,s)$ on $X$ is a coherent sheaf $F$ on $X$
 and a  section $s\in H^0(X,F)$ satisfying the following stability conditions:
\begin{itemize}
\item $F$ is \emph{pure} of dimension 1,
\item the section $s:\OO_X\to F$ has cokernel of dimensional 0.
\end{itemize}
\end{definition}

Let $C$ be the scheme-theoretic support of $F$.
By the purity condition, all the irreducible components
of $C$ are of dimension 1 (no 0-dimensional components are permitted).
By  \cite[Lemma 1.6]{pt}, the kernel of $s$ is the ideal sheaf of $C$,
$$\II_C=\text{ker}(s) \subset \OO_X\, ,$$
and $C$ has no embedded points.
A stable pair $$\OO_X\to F$$ 
therefore
defines
 a Cohen-Macaulay subcurve $C\subset X$ via the kernel of $s$
 and a 0-dimensional subscheme\footnote{When $C$ is Gorenstein (for instance if $C$ lies
in a nonsingular surface), stable pairs supported on $C$ are in bijection with 0-dimensional subschemes of $C$. More precise scheme theoretic isomorphisms of moduli spaces are proved in \cite[Appendix B]{pt3}.}
of $C$ via the support of the
 cokernel of $s$.

To a stable pair, we associate the Euler characteristic and
the class of the support $C$ of $F$,
$$\chi(F)=n\in \mathbb{Z} \  \ \ \text{and} \ \ \ [C]=\beta\in H_2(X,\mathbb{Z})\,.$$
For fixed $n$ and $\beta$,
there is a projective moduli space of stable pairs $P_n(X,\beta)$. 
Unless $\beta$ is an effective curve class, the moduli space
$P_n(X,\beta)$ is empty.

A foundational
treatment of the moduli space
of stable pairs 
is presented in \cite{pt} via the results of  Le Potier \cite{LePot}.
Just as the Hilbert scheme $I_n(X,\beta)$
of subcurves of $X$ of Euler characteristic $n$ and class
$\beta$ is a fine moduli space with a universal quotient sequence, 
$P_n(X,\beta)$ is a fine moduli space with a universal stable pair \cite[Section 2.3]{pt}. While the Hilbert scheme $I_n(X,\beta)$
is a moduli space of curves with free and embedded points, the
moduli space of stable pairs $P_n(X,\beta)$
 should be viewed as a moduli space of curves with points 
\emph{on the curve} determined by the cokernel of $s$. 
Though the additional points
still play a role,
$P_n(X,\beta)$
much smaller than $I_n(X,\beta)$. 

If $P_n(X,\beta)$ is non-empty, then $P_{m}(X,\beta)$ is
non-empty for all $m>n$. Stable pairs with higher Euler characteristic
can be obtained by suitably twisting stable pairs with lower Euler
characteristic (in other words, by {\em adding points}).
On the other hand, for a fixed class $\beta\in H_2(X,\mathbb{Z})$,
the moduli space $P_n(X,\beta)$ is
empty for all sufficiently negative $n$. The proof 
exactly parallels the same result for the Hilbert scheme of
curves $I_n(X,\beta)$.

\subsection{Action of the descendents}\label{actact}
Denote the universal stable pair over $X\times P_{n}(X,\beta)$ by
$$\OO_{X\times P_n(X,\beta)} \stackrel{s\ }{\rightarrow} \FF .$$
For a stable pair $(F,s)\in P_{n}(X,\beta)$, the restriction of
the universal stable pair
to the fiber
 $$X \times (F,s) \ \subset\  
X\times P_{n}(X,\beta)
$$
is canonically isomorphic to $\OO_X\stackrel{s\ }{\rightarrow} F$.
Let
$$\pi_X\colon X\times P_{n}(X,\beta)\to X,$$
$$\pi_P\colon X\times P_{n}(X,\beta)
\to P_{n}(X,\beta)$$
 be the projections onto the first and second factors.
Since $X$ is nonsingular
and
$\FF$ is $\pi_P$-flat, $\FF$ has a finite resolution 
by locally free sheaves.{\footnote{Both $X$ and $P_n(X,\beta)$
carry ample line bundles.}}
Hence, the Chern character of the universal sheaf $\FF$ on 
$X \times P_n(X,\beta)$ is well-defined.

\begin{definition}\label{dact}
For each cohomology{\footnote{All homology
and cohomology groups will be taken with $\Q$-coefficients
unless explicitly denoted otherwise.}}
class $\gamma\in H^*(X)$ and integer $i\in \mathbb{Z}_{\geq 0}$,
the action of the {\em descendent} $\tau_i(\gamma)$ is defined by
$$
\tau_i(\gamma)=\pi_{P*}\big(\pi_X^*(\gamma)\cdot \text{\em ch}_{2+i}(\FF)
\cap \pi_P^*(\ \cdot\ )\big)\, .$$ 
\end{definition}
\noindent The pull-back
$\pi^*_P$ is well-defined in homology since $\pi_p$ is flat \cite{Ast}.

We may view the descendent action as defining a
cohomology class 
$$\tau_i(\gamma)\in H^*(P_n(X,\beta))$$
or as defining an endomorphism
$$\tau_i(\gamma):
H_*(P_{n}(X,\beta))\to H_*(P_{n}(X,\beta))\, .
$$
Definition \ref{dact} is the standard method of obtaining
classes on moduli spaces of sheaves via universal structures.
The construction has been used previously for the
cohomology of the moduli space of bundles on a curve \cite{New},
for the cycle theory of 
the Hilbert schemes of points of a surface \cite{ES}, and in
Donaldson's famous $\mu$ map for gauge theory on 4-manifolds \cite{DonK}.

\subsection{Tautological classes} \label{actactt}
Let $\mathbb{D}$ denote the polynomial $\Q$-algebra on the symbols 
$$\{ \, \tau_i(\gamma)\, |\, i\in \mathbb{Z}_{\geq 0} \text{ and } \gamma\in H^*(X)\, \}$$
subject to the basic  linear relations
\begin{eqnarray*}
\tau_i(\lambda\cdot \gamma) & = & \lambda \tau_i(\gamma)\, ,\\
\tau_i(\gamma+ \widehat{\gamma}) &=& \tau_i(\gamma)+ \tau_i(\widehat{\gamma})\, ,
\end{eqnarray*}
for $\lambda \in \Q$ and $\gamma, \widehat{\gamma}\in H^*(X)$.
The descendent action defines a $\Q$-algebra homomorphism
$$\alpha^X_{n,\beta}: \mathbb{D} \rightarrow H^*(P_n(X,\beta))\, .$$
The most basic questions about the descendent action are
to determine
$$\text{Ker}(\alpha^X_{n,\beta}) \subset \mathbb{D} \ \ \text { and } \ \
\text{Im}(\alpha^X_{n,\beta}) \subset H^*(P_n(X,\beta))\, .$$ 
Both questions are rather difficult since the space $P_n(X,\beta)$
can be very complicated (with serious singularities and
components of different dimensions). Few methods are
available to study $H^*(P_n(X,\beta))$.
 
Following the study of the cohomology of the moduli of stable curves, we 
define, for the moduli space of stable pairs $P_n(X,\beta)$,
\begin{enumerate} 
\item[$\bullet$]$ \text{Im}(\alpha^X_{n,\beta})\subset H^*(P_n(X,\beta))$ to be the algebra of
{\em tautological classes},
\item[$\bullet$] 
$\text{Ker}(\alpha^X_{n,\beta})\subset \mathbb{D}$ to be the  the ideal of {\em tautological relations}.
\end{enumerate}
The basic expectation is that natural constructions yield tautological
classes. For the moduli spaces of curves there is a long
history of the study of tautological classes,
geometric constructions, and relations, see \cite{FaPa, PaSLC} for surveys.

As a simple example, consider the tautological classes in the case 
$$X=\PP^3\, , \ \ \ n=1\, , \ \ \ \beta=\mathsf{L}\, ,$$ 
where $\mathsf{L}\in H_2(\mathsf{P}^3,\mathbb{Z})$ is the class
of a line. The moduli space $P_1(\PP^3,\mathsf{L})$ is isomorphic to the Grassmannian
$\mathbb{G}(2,4)$. The ring homomorphism 
$$\alpha_{1,\mathsf{L}}^{\PP^3}: \mathbb{D} \rightarrow H^*(P_1(\PP^3,
\mathsf{L}))$$
is surjective, so {\em all} classes are tautological.
The tautological relations
 $$\text{Ker}(\alpha_{1,\mathsf{L}}^{\PP^3})\subset \mathbb{D}$$ can  be
determined by the Schubert calculus.

Our study of descendents here follows a different line which 
is more accessible than the full analysis of $\alpha_{n,\beta}^X$.
The moduli space $P_n(X,\beta)$ carries a virtual fundamental
class 
$$[P_n(X,\beta)]^{vir} \in H_*(P_n(X,\beta))$$
obtained from the deformation theory of stable pairs.
There is an associated integration map
\begin{equation}\label{fred}
\int_{[P_n(X,\beta)]^{vir}} :\ \mathbb{D} \rightarrow \Q
\end{equation}
defined by 
$$\int_{[P_n(X,\beta)]^{vir}} \mathsf{D} = \int_{P_n(X,\beta)} \alpha_{n,\beta}^X (\mathsf{D}) \cap
[P_n(X,\beta)]^{vir}\, $$
for $\mathsf{D}\in \mathbb{D}$. Here,
$$\int_{P_n(X,\beta)}:\ H_*(P_n(X,\beta)) \rightarrow \Q$$
is the canonical point counting map factoring through $H_0(P_n(X,\beta))$.
The standard theory of descendents is a study of the
integration map \eqref{fred}.

\subsection{Deformation theory}
To define a virtual fundamental class \cite{BehFan,LiTian}, 
a 2-term deformation/obstruction theory must be found
on the moduli space of stable pairs $P_n(X,\beta)$.
As in the case of the Hilbert scheme
$I_n(X,\beta)$, the most immediate
 obstruction theory of $P_n(X,\beta)$ does \textit{not} admit 
such a structure. For $I_n(X,\beta)$, a suitable
obstruction theory
is obtained by viewing $C\subset X$ {\em not} as
a subscheme, but rather as an ideal 
sheaf $\II_C$ with trivial determinant \cite{DT,THC}. 
For $P_n(X,\beta)$, a suitable obstruction theory is obtained
by viewing
a stable pair {\em not} as sheaf with a section, but
as an object
$$[\OO_X\rightarrow F]\in D^b(X)$$ 
in the bounded derived category of coherent sheaves on $X$.

Denote the quasi-isomorphism equivalence class of
the complex
$[\OO_X\rightarrow F]$ in $D^b(X)$ by
$I\udot$.
The quasi-isomorphism 
class $I\udot$ determines{\footnote{The claims require
the dimension of $X$ to be 3.}}
the stable pair 
\cite[Proposition 1.21]{pt}, and the fixed-determinant deformations 
of $I\udot$ in $D^b(X)$ match those of the pair $(F,s)$ to all orders
\cite[Theorem 2.7]{pt}. The latter property shows 
the scheme $P_n(X,\beta)$ may be viewed as a moduli space
of objects in the derived 
category.{\footnote{The moduli of objects
in the derived category usually yields
Artin stacks. The space $P_n(X,\beta)$ is a
rare example where a component of the moduli of objects
in the derived category is a scheme (uniformly for all $3$-folds $X$).}}
We can then use the obstruction theory of the complex
$I\udot$ rather than the
obstruction theory of sheaves with sections.

The deformation/obstruction theory
for complexes  at  $[I\udot]\in P_n(X,\beta)$ is governed by
\begin{equation}\label{exts2}
\Ext^1(I\udot,I\udot)_0 \quad\mathrm{and}\quad \Ext^2(I\udot,I\udot)_0\,.
\end{equation}
The obstruction theory \eqref{exts2} has all the formal properties
of the Hilbert scheme case: 2 terms, 
a virtual class of (complex) dimension $d_\beta=\int_\beta c_1(X)$, 
$$[P_{n}(X,\beta)]^{vir} \in H_{2d_\beta} 
\big(P_n(X,\beta),\mathbb{Z}\big)\, ,$$
and a description 
via the  $\chi^B$-weighted Euler characteristics in the Calabi-Yau 
case \cite{Kai}.

\subsection{Descendent invariants}\label{dess}
Let $X$ be a nonsingular projective 3-fold.
For nonzero $\beta\in H_2(X,\Z)$ and arbitrary $\gamma_i\in H^*(X)$,
define the stable pairs invariant with descendent insertions by
\begin{equation}\label{ddd}
\Big\langle \tau_{k_1}(\gamma_1)\ldots \tau_{k_r}(\gamma_r)
\Big\rangle_{\!n,\beta}^X  = 
\int_{[P_{n}(X,\beta)]^{vir}}
\prod_{i=1}^r \tau_{k_i}(\gamma_i)\ . 
\end{equation}
The partition function is 
\begin{equation}\label{ppp}
\ZZ_{\mathsf{P}}\Big(X;q\ \Big|   \prod_{i=1}^r \tau_{k_i}(\gamma_{i})
\Big)_\beta
=\sum_{n\in \mathbb{Z}} 
\Big\langle \prod_{i=1}^r \tau_{k_i}(\gamma_{i}) 
\Big\rangle_{\!n,\beta}^X q^n.
\end{equation}
Since $P_n(X,\beta)$ is empty for sufficiently negative
$n$, the partition function 
is a Laurent series in $q$,
$$\bZ_{\mathsf{P}}\Big(X;q\ \Big|\ \prod_{i=1}^r \tau_{k_i}(\gamma_{i})\Big)_\beta \in \mathbb{Q}((q))\, .$$ 

The descendent invariants \eqref{ddd} and the associated partition functions \eqref{ppp}
are the central topics of the paper. From the point of view of the complete tautological
ring of descendent classes on $P_n(X,\beta)$, the descendent invariants \eqref{ddd} 
constitute only small part of the full data. However, among many advantages,
the integrals \eqref{ddd} are deformation invariant as $X$ varies in families.
The same can not be said of the tautological ring nor of the full cohomology
$H^*(P_n(X,\beta))$.

In addition to carrying data about the tautological classes on $P_n(X,\beta)$,
the descendent series are related to the enumerative geometry of curves in $X$.
The connection is clearest for the primary fields $\tau_0(\gamma)$ which correspond
to incidence conditions for the support curve of the stable pair with a fixed cycle
$$V_\gamma \subset X$$
of class $\gamma \in H^*(X)$. But even for primary fields, the partition
function 
$$\bZ_{\mathsf{P}}\Big(X;q\ \Big|\ \prod_{i=1}^r \tau_{0}(\gamma_{i})\Big)_\beta $$
provides a virtual count and is rarely strictly enumerative.

Descendents $\tau_k(D)$, for $k\geq 0$ and $D\subset X$ a divisor, can be
viewed as imposing tangency conditions
of the support curve of the stable pair along the divisor $D$.
The connection of $\tau_k(D)$ to tangency conditions is not as close as the enumerative interpretation
of primary fields --- the tangency condition is just the leading term in the
understanding of $\tau_k(D)$. 
The topic will be discussed in Section \ref{ggogg}.
 
\subsection{Plan of the paper}
The paper starts in Section \ref{111} with a discussion
of the rationality of the descendent partition function
in absolute, equivariant, and relative geometries.
While the general statement is conjectural, rationality in toric and
hypersurface geometries has been proven
in joint work with A. Pixton in \cite{part1,PP2,PPQ}. 
Examples of exact calculations of descendents are given in Section 
\ref{fex}.
A precise conjecture for a functional equation related to
the change of variable $$q \mapsto \frac{1}{q}$$
is presented
in Section \ref{funk}, and a conjecture constraining the
poles appears in Section \ref{polec}.

The second topic, the Gromov-Witten/Pairs 
correspondence for descendents, is discussed in Section \ref{222r}. 
The descendent theory of stable maps
and stable pairs on a nonsingular projective 3-fold $X$
are conjectured to be {\em equivalent} via a universal
transformation. While the correspondence is
proven in joint work with A. Pixton in
toric \cite{PPDC} and hypersurface \cite{PPQ}  cases and
several formal properties are established, a closed 
formula for the transformation is not known.
 
The Gromov-Witten/Pairs correspondence has motivated
much of the development of the descendent theory on
the sheaf side. The first such conjectures for descendent series
 were made in joint work with D. Maulik, A. Okounkov, and N. Nekrasov
\cite{MNOP1,MNOP2} in the context of the 
Gromov-Witten/Donaldson-Thomas correspondence
for the partition functions associated
to the Hilbert schemes $I_n(X,\beta)$ of subcurves of $X$.

Given the Gromov-Witten/Pairs correspondence and the
well-known Virasoro constraints for descendents
in Gromov-Witten theory, there must be corresponding
Virasoro constraints for the descendent theory
of stable pairs. For 
the Hilbert schemes $I_n(X,\beta)$ of curves,
descendent constraints were studied by A. Oblomkov, A. Okounkov,
and myself in Princeton a decade ago \cite{oop}.
In Section \ref{333r},  conjectural descendent constraints
for the stable pairs theory of $\PP^3$ are presented (joint
work with A. Oblomkov and A. Okounkov).

The moduli space of stable pairs $P_n(X,\beta)$
has a virtual fundamental class in homology $H_*(P_n(X,\beta))$.
By construction, the class lifts to algebraic cycles $A_*(P_n(X,\beta))$.
In a recent paper, Junliang Shen has lifted the virtual fundamental
class further to algebraic cobordism $\Omega_*(P_n(X,\beta))$.
Shen's results open a new area of exploration with beautiful
structure. At the moment, the methods available to explore
the virtual fundamental class in cobordism all use the theory of descendents
(since the Chern classes of the virtual tangent bundle of
$P_n(X,\beta)$ are {\em tautological}).
Shen's work is presented in Section \ref{4448}.

\subsection{Acknowledgments}
Discussions with J. Bryan, S. Katz, D. Maulik, G. Oberdieck,
A. Oblomkov, A. Okounkov,
A. Pixton, J. Shen, 
R. Thomas, Y. Toda, and Q. Yin 
about  stable pairs
and descendent invariants
have played an important role in my view of the subject.
I was partially supported by
 SNF grant 200021-143274,   ERC grant
AdG-320368-MCSK, SwissMAP, and the Einstein Stiftung.

The perspective of the paper is based in part on my talk {\em Why descendents?}
at the Newton institute in Cambridge in the spring of 2011, though
much of the progress discussed here has happened since then.

\section{Rationality} \label{111r}
\subsection{Overview}
Let $X$ be a nonsingular projective 3-fold. Our goal here is to present the
conjectures governing the {\em rationality} of the partition functions of
descendent invariants for the stable pairs theory of $X$.
The most straightforward statements are for the absolute theory, but we
will present the rationality claims for the equivariant
and relative stable pairs theories as well. The latter two 
appear naturally when studying the absolute theory: most results to date
involve equivariant and relative techniques.
In addition to rationality, we will also discuss the 
{\em functional equation} and the {\em pole constraints} for the
descendent partition functions. 

While rationality has been established in many cases, new ideas 
are required to prove the conjectures in full generality. The subject
intertwines the Chern characters of the universal sheaves with the
geometry of the virtual fundamental class. Perhaps, in the future,
a point of view will emerge from which rationality is obvious.
Hopefully, the 
 functional equation will then also be clear. At present, 
the geometries for which the functional equation has been
proven are rather few. 

\subsection{Absolute theory}
Let $X$ be a nonsingular projective 3-fold. The stable pairs theory for $X$ as
presented in the introduction is the {\em absolute} case.
Let $\beta\in H_2(X,\Z)$ be a nonzero class, and let  $\gamma_i\in H^*(X)$.
The following conjecture{\footnote{A weaker conjecture
for descendent partition functions for the Hilbert scheme $I_n(X,\beta)$
was proposed earlier in \cite{MNOP2}.}}  was proposed{\footnote{Theorems and Conjectures 
are dated in the text by the year of the arXiv posting. The published
dates are later and can be found in the bibliography.}} in \cite{pt2}.
 
\begin{Conjecture}[P.-Thomas, 2007]
\label{111} 
For $X$ a nonsingular projective $3$-fold,
the descendent partition function
$$\ZZ_{\mathsf{P}}\big(X;q\ |   \prod_{i=1}^r \tau_{k_i}(\gamma_{i})
\big)_\beta$$ is the 
Laurent expansion in $q$ of a rational function in $\mathbb{Q}(q)$.
\end{Conjecture}

In the absolute case, the descendent series satisfies a dimension constraint.
For $\gamma_i\in H^{e_i}(X)$, the (complex) degree of the insertion 
$\tau_{k_i}(\gamma_i)$ is $\frac{e_i}{2}+k_i-1$.
If the sum of the degrees of the descendent insertions does
 not equal the virtual dimension,
$$\text{dim}_\com\, [ P_n(X,\beta)]^{vir} = \int_\beta c_1(X)\,, $$
the partition function $\ZZ_{\mathsf{P}}\big(X;q\ |   \prod_{i=1}^r \tau_{k_i}(\gamma_{i})
\big)_\beta$ vanishes.

In case $X$ is a nonsingular projective Calabi-Yau 3-fold, the virtual dimension
of $P_n(X,\beta)$ is always 0. The rationality of the basic partition function
         $$\ZZ_{\mathsf{P}}\big(X;q\ |\,   \mathsf{1} \big)_\beta$$
was proven{\footnote{See \cite{pt3}
for a similar rationality argument in a restricted (simpler) setting.}} 
in \cite{Bridge,Toda} by Serre duality, wall-crossing, and a weighted
Euler characteristic approach to the virtual class \cite{Kai}. 
At the moment, the proof for Calabi-Yau 3-folds does not appear to
suggest an approach in the general case.

\subsection{Equivariant theory}
Let $X$ be a nonsingular quasi-projective toric 3-fold equipped
with an action of the 3-dimensional torus
$$\T= \com^* \times \com^* \times \com^*\, .$$ 
The stable pairs descendent
invariants can be lifted to equivariant cohomology (and
defined by residues in the open case). For
equivariant classes $\gamma_i \in H^*_{\T}(X)$, the
descendent partition function is a Laurent series in $q$,
$$\bZ_{\mathsf{P}}\Big(X;q\ \Big|\ \prod_{i=1}^r \tau_{k_i}(\gamma_{i})\Big)^\T_\beta \in 
\mathbb{Q}(s_1,s_2,s_3)((q))\, ,$$ 
with coefficients in the field of fractions of  
$$H^*_\T(\bullet)=\mathbb{Q}[s_1,s_2,s_3]\,.$$
The stable pair theory for such toric $X$ is the {\em equivariant} case.
A central result of \cite{part1,PP2} is the following rationality
property.

\begin{Theorem}[P.-Pixton, 2012] \label{pp12}
For $X$ a nonsingular quasi-projective toric 3-fold,
 the descendent partition function
$$\ZZ_{P}\Big(X;q\ \Big|   \prod_{i=1}^r \tau_{k_i}(\gamma_i)
\Big)^\T_\beta$$ is the 
Laurent expansion in $q$ of a rational function in 
$\mathbb{Q}(q,s_1,s_2,s_3)$.
\end{Theorem}

The proof of Theorem \ref{pp12} uses the virtual localization
formula of \cite{GraberP}, the capped vertex{\footnote{A basic tool in the proof is the capped {\em descendent} vertex.
The 1-leg capped descendent vertex is proven to be rational in \cite{part1}.
The 2-leg and 3-leg capped descendent vertices are proven to be rational
in \cite{PP2}.}}  
perspective of
\cite{moop}, 
the quantum cohomology of the 
Hilbert scheme of points of resolutions of $A_r$-singularities \cite{mo1,hilb1}, 
and a delicate argument for pole cancellation at the vertex
\cite{part1}.
In the toric case, calculations can be made effectively,
but the computational methods are not very efficient.

When $X$ is a nonsingular projective
toric $3$-fold, Theorem  \ref{pp12}
implies Conjecture \ref{111} for $X$ by taking the non-equivariant limit. However, Theorem \ref{pp12} is much 
stronger in the toric case than Conjecture \ref{111} 
since the descendent insertions may exceed the virtual dimension in equivariant cohomology.

In addition to the Calabi-Yau and toric cases, Conjecture \ref{111}
has been proven in \cite{PPQ} for complete intersections in products of 
projective spaces (for descendents of cohomology classes $\gamma_i$ 
restricted from the ambient space --- the precise statement
is presented in  Section \ref{compint}).  
Taken together, the evidence for Conjecture \ref{111} is
rather compelling.

\subsection{First examples}\label{fex}
Let $X$ be a nonsingular projective Calabi-Yau 3-fold, and
let
$$C \subset X\, $$
be a rigid nonsingular {\em rational} curve.
Let 
$$\ZZ_{\mathsf{P}}\big(C\subset X;q\ |\,   \mathsf{1} \big)_{d[C]}$$
be the contribution to the partition function
$\ZZ_{\mathsf{P}}\big(X;q\ | \, \mathsf{1} \big)_{d[C]}$
obtained from the moduli of stable pairs {\em supported on} $C$.
A localization calculation which goes back to the Gromov-Witten
evaluation of \cite{FP} yields
\begin{equation}\label{jqq}
\ZZ_{\mathsf{P}}\big(C\subset X;q\ |   1 \big)_{d[C]}=\sum_{\mu\vdash d}
\frac{(-1)^{\ell(\mu)}}{{\mathfrak{z}}(\mu)} \prod_{i=1}^{\ell(\mu)} \frac{(-q)^{m_i}}
{(1-(-q)^{m_i})^2}\, .
\end{equation}
The sum here is over all (unordered) partitions of $d$,
$$\mu=(m_1,\ldots,m_{\ell(\mu)})\, , \ \ \ \sum_{i=1}^{\ell(\mu)}m_i = d\, ,$$
and ${\mathfrak{z}}(\mu)$ is the standard combinatorial factor
$${\mathfrak{z}}(\mu)= \prod_{i=1}^{\ell(\mu)} m_i \cdot |\text{Aut}(\mu)|\, .$$ 
The evaluation \eqref{jqq} played an important role
in the discovery of the Gromov-Witten/Donaldson-Thomas
correspondence in \cite{MNOP1}.

In example \eqref{jqq}, only the trivial descendent insertion $\mathsf{1}$ 
appears. 
For non-trivial insertions, consider the case where $X$ is $\PP^3$. Let
$$\mathsf{p},\mathsf{L}\in  H_*(\PP^3)$$ 
be the point and line classes  in $\PP^3$ respectively.
Geometrically, there is unique line through two points of $\PP^3$.
The corresponding partition function is also simple,
\begin{equation}\label{ggtt}
\ZZ_{\mathsf{P}}\big(\PP^3;q\ |\,   \tau_0(\mathsf{p})\tau_0(\mathsf{p}) \big
)_{\mathsf{L}} =  q+2q^2+q^3\, .
\end{equation}
The resulting series is not only rational, but in fact 
polynomial. 
For curve class $\mathsf{L}$,
the descendent invariants in \eqref{ggtt} vanish
for Euler characteristic greater than 3.

In example \eqref{ggtt}, only primary fields  (with descendent subscript 0) 
appear.
An example with higher descendents is
$$\ZZ_{\mathsf{P}}\big(\PP^3;q\ |\,   \tau_2(\mathsf{p}) \big
)_{\mathsf{L}} =  \frac{1}{12}q-\frac{5}{6}q^2+\frac{1}{12}q^3\, .$$
The fractions here come from the Chern character. Again, the
result is a cubic polynomial. 
More interesting is the partition function
\begin{equation}\label{s555}
\ZZ_{\mathsf{P}}\big(\PP^3;q\ |\,   \tau_5(\mathsf{1}) \big
)_{\mathsf{L}} =  \frac{-2q-q^2+31q^3-31q^4+q^5+2q^6}{18(1+q)^3}\, .
\end{equation}

The partition functions considered so far are all in
the absolute case. For an equivariant descendent series,
consider the $\T$-action on $\PP^3$ defined
by representation
weights $\lambda_0,\lambda_1,\lambda_2,\lambda_3$
on the vector space $\com^4$.
Let $$\mathsf{p_0}\in H^4_\T(\PP^3)$$
be the class of the $\T$-fixed point corresponding
to the weight $\lambda_0$ subspace of $\com^4$.
Then, 
$$\ZZ_{\mathsf{P}}\big(\PP^3;q\ |   \tau_3(\mathsf{p}_0) \big
)_{\mathsf{L}} = 
\frac{\mathsf{A} q-\mathsf{B} q^2+\mathsf{B}q^3-\mathsf{A}q^4}{(1+q)}
\, $$
where $A,B\in H^2_\T(\bullet)$ are given by
\begin{eqnarray*}
\mathsf{A}& =& \frac{1}{8}\lambda_0 - \frac{1}{24}(\lambda_1+\lambda_2+\lambda_3)\, , \\
\mathsf{B} & = &  \frac{9}{8}\lambda_0 -\frac{3}{8} (\lambda_1+\lambda_2+\lambda_3)\, .
\end{eqnarray*}
The descendent insertion here has dimension 5 which
exceeds the virtual dimension 4 of the moduli space
of stable pair, so the invariants lie in $H^2_\T(\bullet)$.
The obvious symmetry in all of these descendent series
is explained by the conjectural function equation (discussed
in Section \ref{funk}).

All of the formulas discussed above are calculated by
the virtual localization formula \cite{GraberP} for stable pairs.
The $\T$-fixed points, virtual tangent weights, and virtual 
normal weights are described in detail in \cite{pt2}.

\subsection{Example in degree 2} \label{dd22}
A further example in the absolute case is the degree 2 series 
{\footnotesize{$\mathsf{Z}_{\mathsf{P}}\big(\mathsf{P}^3; q \ |\, \tau_9(1)\big)_{2\mathsf{L}}$}}.
While a rigorous answer could be obtained, the available
computer calculation
here  outputs a conjecture,{\footnote{The answer relies on an old  program 
for the theory of ideal sheaves written by A. Okounkov and
a newer DT/PT descendent correspondence \cite{oop}.}}

{\footnotesize{
$$\hspace{-280pt}\mathsf{Z}_{\mathsf{P}}\big(\mathsf{P}^3; q \ |\, \tau_9(1)\big)_{2\mathsf{L}}
=$$}}
{\tiny{ \vspace{5pt}
$$
-\frac{(73 q^{12}   - 825 q^{11}   - 124 q^{10}   + 5945 q^{9}  + 779 q^{8}  - 36020 q^{7}  + 60224 q^{6}  - 36020 q^{5}  + 779 q^{4}  + 5945 q^{3}  - 124 q^2  - 825 q
+ 73) q}    {60480 (1 + q)^3  (-1 + q)^3 }\, .$$}}
             
\noindent The computer calculations of Section \ref{fex} all provide rigorous
results and could be improved to handle higher degree curves, but the code
has not yet been written. 

\subsection{Relative theory}\label{relth}
Let $X$ be a nonsingular projective $3$-fold containing
a nonsingular divisor
$$D\subset X\, .$$ 
The {\em relative} case concerns the geometry $X/D$. 
%

The moduli space 
$P_{n}(X/D,\beta)$
parameterizes
stable relative pairs
\begin{equation}\label{vyq}
s:\OO_{X[k]} \rightarrow F
\end{equation}
on the $k$-step{\footnote{We
follow the terminology of \cite{L, LW}.}} degeneration $X[k]$.

\vspace{8pt}
\noindent $\bullet$  The algebraic variety $X[k]$ is constructed
by attaching a chain of $k$ copies
of the 3-fold
$\mathbb{P}(N_{X/D}\oplus \OO_D)$ 
equipped with 0-sections and $\infty$-sections
$$D \stackrel{\iota_0} \longrightarrow \mathbb{P}(N_{X/D}\oplus \OO_D)
\stackrel{\ \iota_\infty}{\longleftarrow} D$$
defined by the summands $N_{X/D}$ and $\OO_D$ respectively.
The $k$-step degeneration $X[k]$ is a union
$$X \cup_D \mathbb{P}(N_{X/D}\oplus \OO_D)\cup_D 
\mathbb{P}(N_{X/D}\oplus \OO_D) \cup_D \cdots \cup_D 
\mathbb{P}(N_{X/D}\oplus \OO_D)\, ,$$
where the attachments are made along $\infty$-sections on
the left and $0$-sections on the right.  
The original divisor $D\subset X$ is considered
an $\infty$-section for the attachment rules.
The rightmost component of $X[k]$ carries
the last $\infty$-section, $$D_\infty \subset X[k],$$
called the {\em relative divisor}.
The $k$-step degeneration also admits a canonical contraction
map
\begin{equation}\label{vssv}
X[k] \rightarrow X
\end{equation}
collapsing all the attached components to $D\subset X$.

\vspace{8pt}
\noindent $\bullet$ The sheaf $F$ on $X[k]$ is of Euler
characteristic
$$\chi(F)=n$$
and has 1-dimensional support on $X[k]$ which pushes-down
via the contraction \eqref{vssv} 
to the class
$$\beta\in H_2(X,\Z).$$

\vspace{8pt}
\noindent $\bullet$ The following stability conditions
are required for stable relative pairs:
\begin{enumerate}
\item[(i)] $F$ is pure with finite locally free resolution,
\item[(ii)] the higher derived functors of the
restriction of $F$ to the singular{\footnote{The singular
loci of $X[k]$ , by convention, include also the
relative divisor $D_\infty\subset X[k]$ even though
$X[k]$ is nonsingular along $D_\infty$ as a variety.
The perspective of log geometry is more natural here.}} loci of $X[k]$ vanish,
\item[(iii)] the section $s$ has 0-dimensional cokernel supported
away from the singular loci of $X[k]$.
\item[(iv)] the pair \eqref{vyq} has only finitely many automorphisms covering
the automorphisms of $X[k]/X$.
\end{enumerate}

\vspace{8pt}
The moduli space $P_n(X/D,\beta)$ of stable relative pairs
is a complete Deligne-Mumford stack
equipped with a map to the Hilbert scheme of points of $D$
via the restriction of the pair to the relative divisor,
$$P_n(X/D,\beta) \to \Hilb(D,\int_\beta [D])\ .$$ 
Cohomology classes on $\Hilb(D,\int_\beta [D])$ may 
thus
be pulled-back to the moduli space $P_n(X/D,\beta)$.

We will use
the \emph{Nakajima basis} of $H^*(\Hilb(D,\int_\beta [D]))$ indexed
by a partition $\mu$ of $\int_\beta [D]$
labeled by cohomology classes of $D$. For example, the
class
$$
\left.\big|\mu\rang \in H^*(\Hilb(D,\int_\beta [D]))\,,
$$
with all cohomology labels equal to the identity,
 is $\prod \mu_i^{-1}$ times
the Poincar\'e dual of the closure of the subvariety formed by unions of
schemes of length
$$
\mu_1,\dots, \mu_{\ell(\mu)}
$$
supported at $\ell(\mu)$ distinct points of $D$.

The stable pairs descendent
invariants in the relative case
are defined using the universal sheaf
just as in the absolute case.
The universal sheaf is defined here on the
universal degeneration of $X/D$ over
$P_n(X/D,\beta)$. 
The cohomology classes $\gamma_i\in H^*(X)$
are pulled-back to the universal degeneration
via the contraction map \eqref{vssv}.
The descendent partition function with boundary conditions
$\mu$ is a Laurent series in $q$,
$$\bZ_{\mathsf{P}}
\Big( X/D;q\ \Big|\, \ \prod_{i=1}^r \tau_{k_i}(\gamma_{i})
\, \Big|\, \mu \Big)_\beta \in \mathbb{Q}((q))\, .$$
The basic rationality statement here is parallel to the 
absolute and equivariant cases.

\begin{Conjecture}
\label{222} 
For $X/D$ a nonsingular projective relative 3-fold,
the descendent partition function
$$\bZ_{\mathsf{P}}
\Big( X/D;q\ \Big|\, \prod_{i=1}^r \tau_{k_i}(\gamma_{i})
\, \Big|\, \mu \Big)_\beta \in \mathbb{Q}((q))\, $$
is the 
Laurent expansion in $q$ of a rational function in $\mathbb{Q}(q)$.
\end{Conjecture}

In case $X$ is a nonsingular quasi-projective toric 3-fold
and $D\subset X$ is a toric divisor, an {\em equivariant
relative} stable pairs theory can be defined.
The rationality conjecture then takes the form expected
by combining the rationality statements in 
the equivariant and relative cases.

\begin{Conjecture}
\label{333} 
For $X/D$ a nonsingular quasi-projective relative toric 3-fold,
the descendent partition function
$$\bZ_{\mathsf{P}}
\Big( X/D;q\ \Big|\, \prod_{i=1}^r \tau_{k_i}(\gamma_{i})
\, \Big|\, \mu \Big)^\T_\beta \in \mathbb{Q}(s_1,s_2,s_3)(q)\, $$
is the 
Laurent expansion in $q$ of a rational function in $\mathbb{Q}(q,s_1,
s_2,s_3)$. 
\end{Conjecture}

Of course, both $\gamma_i\in H^\bullet_\T(X)$
and the Nakajima basis element $$\mu \in H^*_\T(\Hilb(D,\int_\beta [D]))$$
must be taken here in equivariant cohomology.
While the full statement of Conjecture \ref{333} remains open,
a partial result follows from Theorem \ref{pp12} and \cite[Theorem 2]{part1}
which addresses the non-equivariant limit in the projective 
relative
toric case.

\begin{Theorem}[P.-Pixton, 2012] \label{ppp12}
For $X/D$ a nonsingular projective relative toric 3-fold,
 the descendent partition function
$$\ZZ_{\mathsf{P}}\Big(X;q\ \Big|   \prod_{i=1}^r \tau_{k_i}(\gamma_i)
\, \Big|\, \mu \Big)_\beta$$ is the 
Laurent expansion in $q$ of a rational function in 
$\mathbb{Q}(q)$.
\end{Theorem}

As an example of a computation in closed form in
the equivariant relative case, consider
the geometry of the {\em cap},
$$\com^2 \times \mathsf{P}^1 / \com^2_\infty\, ,$$
where $\com^2_\infty \subset \com^2 \times \mathsf{P}^1$
is the fiber of 
 $$\com^2 \times \mathsf{P}^1 \rightarrow \mathsf{P}^1$$
over $\infty \in \mathsf{P}^1$.
The first two factors of the 3-torus $\T$ act on the
$\com^2$-factor of the cap with tangent weights $-s_1$ and $-s_2$.
The third factor of $\T$ acts on $\PP^1$ factor of the cap with tangent weights
$-s_3$ and $s_3$ at $0\in \PP^1$ and $\infty\in \PP^1$
respectively.

From several perspectives, the equivariant
relative descendent partition function 
$$\mathsf{Z}^{\mathsf{cap}}_{\mathsf{P}}( \, \tau_d(\mathsf{p})\, |\,(d) )^\T_d =
\sum_{n}
\Big\langle \tau_{d}(\mathsf{p})\, \Big| \, (d) 
\Big\rangle_{\!n,d}^{\! \text{Cap}}\, q^n\ , \ \ \ \  d>0 \ $$
is the most important series in the cap geometry \cite{PPstat}.
Here, $$\mathsf{p}\in H_\T^2(\com^2 \times \PP^1)$$ is the
class of the $\T$-fixed point of
$\com^2\times \PP^1$ over $0\in \PP^1$, and
the Nakaijima basis element $(d)$ is weighted with the
identity class in $H^*_\T(\text{Hilb}(\com^2,d))$.
A central result of \cite{PPstat} is the following calculation.{\footnote{The
formula here differs from \cite{PPstat} by a factor of $s_1s_2$
since a different convention for the cohomology class $\mathsf{p}$
is taken.}}

\begin{Theorem}[P.-Pixton, 2011] \label{yty7}
We have
$$
\mathsf{Z}^{\mathsf{cap}}_{\mathsf{P}}( \, \tau_d(\mathsf{p})\, |\,(d) )^\T_d
 =
\frac{q^d}{d!}\left(\frac{s_1+s_2}{2}\right)
\sum_{i=1}^d  \frac{ 1+(-q)^{i}}{1-(-q)^i} \ . $$
\end{Theorem}

In the above formula, the coefficient of $q^d$, 
$$ \big\langle \tau_d(\mathsf{p}), (d) \big\rangle_{\text{Hilb}(\com^2,d)}=
\frac{s_1+s_2}{2\cdot (d-1)!}\, ,$$
is the classical $(\com^*)^2$-equivariant pairing on the
Hilbert scheme of points $\text{Hilb}(\mathbb{C}^2,d)$.
The proof of Theorem \ref{yty7}
is a rather delicate localization calculation (using several
special properties such as the \`a priori divisibility of the answer by 
$s_1+s_2$ from the holomorphic symplectic form on 
$\text{Hilb}(\mathbb{C}^2,d)$).

The difficulty in Theorem \ref{yty7} is obtaining a closed
form evaluation for all $d$. Any particular descendent
series can be calculated by the localization methods.
A calculation, for example, {\em not} covered by Theorem \ref{yty7} is
\begin{multline}\label{fvvt}
\mathsf{Z}^{\mathsf{cap}}_{\mathsf{P}}( \, \tau_2(\mathsf{p})\, |\,(1) )^\T_1
 = 
\big({2s_1^2+3s_1s_2+2s_2^2}\big)q\frac{(1+q^2)}{(1+q)^2} \\ 
+ \big({6s_3(s_1+s_2)-2s_1^2-6s_1s_2-2s^2_2}\big)\frac{q^2}{(1+q)^2}\, .
\end{multline}
A simple closed formula for all descendents of the cap
is unlikely to exist.

\subsection{Functional equation}\label{funk}
In case $X$ is a nonsingular Calabi-Yau 3-fold, the descendent
series viewed as a rational function in $q$ satisfies the
symmetry
\begin{equation}\label{symm}
\ZZ_{\mathsf{P}}\big(X;\frac{1}{q}\, \big|\,   1 \big)_\beta
= \ZZ_{\mathsf{P}}\big(X;{q}\, \big|\,   1 \big)_\beta \, 
\end{equation}
as conjectured in \cite{MNOP1,pt} and proven in \cite{Bridge,Toda}.
In fact, a functional equation for the descendent
partition function is expected to hold
in {\em all} cases (absolute, equivariant, and relative).
For the relative case, the functional equation  is given by the following
formula{\footnote{The conjecture is stated in \cite{part1,PPstat} with a 
sign error: the factor of $q^{d_\beta}$ on the right side of
the functional equation \cite{part1,PPstat} should be $(-q)^{d_\beta}$. Then
two factors of $(-1)^{d_\beta}$ multiply to $1$ and yield 
Conjecture \ref{444} as stated here.}}
\cite{part1,PPstat}.

\begin{Conjecture}[P.-Pixton, 2012] \label{444}
For $X/D$ a nonsingular projective relative 3-fold,
 the descendent series viewed as a rational function in $q$
satisfies the functional equation
$$\ZZ_{\mathsf{P}}\Big(X;\frac{1}{q}\ \Big|   \prod_{i=1}^r \tau_{k_i}(\gamma_i)
\, \Big|\, \mu \Big)_\beta
= (-1)^{|\mu|-\ell(\mu)
+
\sum_{i=1}^r
k_i} q^{-d_\beta}   \ZZ_{\mathsf{P}}\Big(X;{q}\ \Big|   \prod_{i=1}^r \tau_{k_i}(\gamma_i) \, \Big|\, \mu \Big)_\beta\, $$
where the constants are 
$$|\mu|=\int_\beta D\,, \ \ \ \ell(\mu)= {\text{\em length}}(\mu)\,, \ \ \ 
d_\beta = \int_{\beta}c_1(X)\, . $$
\end{Conjecture}

The functional equation in the absolute case is obtained
by specializing  the divisor $D\subset X$ to the empty set
in Conjecture \ref{444}:
$$\ZZ_{\mathsf{P}}\Big(X;\frac{1}{q}\ \Big|   \prod_{i=1}^r \tau_{k_i}(\gamma_i)
\Big)_\beta
= (-1)^{
\sum_{i=1}^r
k_i} q^{-d_\beta}\,   \ZZ_{\mathsf{P}}\Big(X;{q}\ \Big|   \prod_{i=1}^r \tau_{k_i}(\gamma_i) \, \Big)_\beta\, . $$
The functional equation in the equivariant case
is conjectured to be identical,
$$\ZZ_{\mathsf{P}}\Big(X;\frac{1}{q}\ \Big|   \prod_{i=1}^r \tau_{k_i}(\gamma_i)
\Big)^\T_\beta
= (-1)^{
\sum_{i=1}^r
k_i} q^{-d_\beta}\,   \ZZ_{\mathsf{P}}\Big(X;{q}\ \Big|   \prod_{i=1}^r \tau_{k_i}(\gamma_i) \, \Big)^\T_\beta\, .$$
Finally, in the equivariant relative case, the functional equation
expected to be same as in Conjecture \ref{444}.

As an example, the descendent series for the cap
evaluated in Theorem \ref{yty7} satisfies the
conjectured functional equation:
\begin{eqnarray*}
\mathsf{Z}^{\mathsf{cap}}_{\mathsf{P}}\left(\frac{1}{q};
 \, \tau_d(\mathsf{p})\, \Big|\,(d) \right)^\T_d
& = &
\frac{q^{-d}}{d!}\left(\frac{s_1+s_2}{s_1s_2}\right)
\frac{1}{2}\sum_{i=1}^d  \frac{ 1+(-q)^{-i}}{1-(-q)^{-i}} \\
& = &
\frac{1}{q^{2d}}\frac{q^{d}}{d!}\left(\frac{s_1+s_2}{s_1s_2}\right)
\frac{1}{2}\sum_{i=1}^d  \frac{ (-q)^{i}+1}{(-q)^{i}-1} \\
& = &
\frac{(-1)^{d-1+d}}{q^{2d}}
\mathsf{Z}^{\mathsf{cap}}_{\mathsf{P}}\big(q;
 \, \tau_d(\mathsf{p})\, |\,(d) \big)^\T_d\, .
\end{eqnarray*}
Here, the constants for the exponent of $(-1)$ in the
functional equation are 
$$|(d)|=d\,, \ \ \ \ell(d)= 1\,, \ \ \ 
d_\beta = 2d\, . $$
It is straightforward to check the functional equation
in all the examples of Section \ref{fex} - \ref{dd22}.

The evidence for the functional equation for
descendent series
 is not as
large as for the rationality.
For the equivariant relative cap, the functional
equation is proven in \cite{PPstat} for all 
descendents series
$$\mathsf{Z}^{\mathsf{cap}}_{\mathsf{P}}\left( \, \prod_{i=1}^r\tau_{k_i}(\mathsf{p})\, \Big|\,(\mu) \right)^\T_d$$
 {\em after} the specialization $s_3=0$.
The
predicted functional equation for
 $$\mathsf{Z}^{\mathsf{cap}}_{\mathsf{P}}( \, \tau_2(\mathsf{p})\, |\,(1) )^\T_1$$ 
 {\em before} the specialization $s_3=0$ 
can be easily checked from the
formula \eqref{fvvt}.
The functional equation is also known
to hold for a special classes of descendent insertions
in the nonsingular projective toric case \cite{PPDC} 
as will be discussed in Section \ref{eee999}.

\subsection{Pole constraints} \label{polec}

Let $X$ be a nonsingular projective 3-fold,
and let $\beta\in H_2(X,\mathbb{Z})$ be an nonzero
class.
 For $\beta$ to be an effective curve class,
the image of $\beta$ in the lattice
\begin{equation}\label{tttt}
H_2(X,\mathbb{Z})/\text{\text{torsion}}
\end{equation}
must also be nonzero. Let $\text{div}(\beta)\in \mathbb{Z}_{>0}$
be the divisibility of the image of $\beta$ in
the lattice \eqref{tttt}.

\begin{Conjecture} \label{555}
For $d=\text{\em div}(\beta)$,
the poles in $q$ of the rational function
$$\ZZ_{P}\Big(X;q\ \Big|   \prod_{i=1}^r \tau_{k_i}(\gamma_i)
\Big)_\beta$$
may occur only at $q=0$ and 
the roots of the polynomials
$$\{ \, 1-(-q)^m \, | \, 1 \leq m \leq d\, \}.$$
\end{Conjecture}

Of the above conjectures, the evidence for
Conjecture \ref{555} is the weakest. 
The prediction is based on a study of the
stable pairs theory of local curves where
the above pole restrictions are always found.
For example,
the evaluation of Theorem \ref{yty7}
is consistent with the pole statement (even though
Theorem \ref{yty7} concerns the  equivariant relative case).
A promotion of Conjecture \ref{555} to cover all
cases also appears reasonable.

\subsection{Complete intersections}\label{compint}
Rationality results for non-toric 3-folds are proven in \cite{PPQ}  by degeneration
methods for several geometries. The simplest
to state concern nonsingular complete intersections of ample divisors
$$ X \subset \PP^{n_1} \times \cdots \times \PP^{n_m}\ .$$

\begin{Theorem}[P.-Pixton, 2012]
\label{qqq111f} 
Let $X$ be a nonsingular Fano or Calabi-Yau complete intersection 3-fold in
a product of projective spaces.
For even classes 
$\gamma_i \in H^{2*}(X)$, the descendent partition function
$$\ZZ_{\mathsf{P}}\Big(X;q\ \Big| \prod_{i=1}^r \tau_{k_i}(\gamma_i)
\Big)_\beta$$
is the Laurent expansion of a rational function in $\mathbb{Q}(q)$.
\end{Theorem}

By the Lefschetz hyperplane result, the even cohomology of such $X$
is exactly the image of the restricted cohomology from the
product of projective spaces.
Theorem \ref{qqq111f} does {\em not} cover the primitive cohomology in
$H^3(X)$. 
Moreover, even for descendents of the even cohomology 
$H^{2*}(X)$ the functional equation and pole conjectures
are open.

\section{Gromov-Witten/Pairs correspondence} \label{222r}

\subsection{Overview}
Let $X$ be a nonsingular projective variety.
Descendent classes
on the moduli spaces of stable maps $\overline{M}_{g,r}(X,\beta)$ in Gromov-Witten theory,
defined using 
cotangent lines at the marked points,
have played a central role
since the beginning of the subject in the early 90s.  Topological
recursion relations,  $J$-functions, and  Virasoro
constraints all essentially concern descendents. 
The importance of descendents in Gromov-Witten theory was hardly a
surprise: cotangent lines on the moduli
spaces  $\overline{M}_{g,r}$ of stable curves were basic to 
their geometric study
 before Gromov-Witten theory
was developed.

In case $X$ is a nonsingular projective {\em 3-fold},
descendent invariants are defined for both  
 Gromov-Witten theory and the theory of stable pairs.
The geometric constructions are rather different, but
a surprising correspondence conjecturally holds: the two descendent
theories are related by a universal correspondence
for {\em all} nonsingular projective
3-folds. In order words, the two descendent theories
contain exactly the same data.

The origin of the Gromov-Witten/Pairs correspondence is found
in the study of ideal sheaves in \cite{MNOP1, MNOP2}.
Since the descendent theory of stable pairs is much better behaved,
the results and conjectures take a better form for stable pairs
\cite{PPDC,PPQ}.

The rationality results and conjectures of Section \ref{111r} are needed
for the statement of the Gromov-Witten/Pairs correspondence. Just
as in Section \ref{111r}, we present the absolute, equivariant, and
relative cases. A more subtle discussion of diagonals is required
for the relative case.

\subsection{Descendents in Gromov-Witten theory}
Let $X$ be a nonsingular projective 3-fold.
Gromov-Witten theory is defined via integration over the moduli
space of stable maps.
Let
 $\overline{M}_{g,r}(X,\beta)$ denote the moduli space of
$r$-pointed stable maps from connected genus $g$ curves to $X$ representing the
class $\beta\in H_2(X, \Z)$. Let 
$$\text{ev}_i: \overline{M}_{g,r}(X,\beta) \rarr X\, ,$$
$$ \LL_i \rarr \overline{M}_{g,r}(X,\beta)$$
denote the evaluation maps and the cotangent line bundles associated to
the marked points.
Let $\gamma_1, \ldots, \gamma_r\in H^*(X)$, and
let $$\psi_i = c_1(\LL_i) \in H^2(\overline{M}_{g,n}(X,\beta))\, .$$
The {\em descendent fields}, denoted by $\tau_k(\gamma)$, correspond 
to the classes $\psi_i^k \text{ev}_i^*(\gamma)$ on the moduli space
of stable maps. 
Let
$$\Big\langle \tau_{k_1}(\gamma_{1}) \cdots
\tau_{k_r}(\gamma_{r})\Big\rangle_{g,\beta} = \int_{[\overline{M}_{g,r}(X,\beta)]^{vir}} 
\prod_{i=1}^r \psi_i^{k_i} \text{ev}_i^*(\gamma_{_i})$$
denote the descendent
Gromov-Witten invariants. Foundational aspects of the theory
are treated, for example, in \cite{BehFan, LiTian}.

Let $C$ be a possibly disconnected curve with at worst nodal singularities.
The genus of $C$ is defined by $1-\chi(\oh_C)$. 
Let $\overline{M}'_{g,r}(X,\beta)$ denote the moduli space of maps
with possibly {disconnected} domain
curves $C$ of genus $g$ with {\em no} collapsed connected components.
The latter condition requires 
 each connected component of $C$ to represent
a nonzero class in $H_2(X,{\mathbb Z})$. In particular, 
$C$ must represent a {nonzero} class $\beta$.

We define the descendent invariants in the disconnected 
case by
$$\Big\langle \tau_{k_1}(\gamma_{1}) \cdots
\tau_{k_r}(\gamma_{r})\Big\rangle'_{g,\beta} = \int_{[\overline{M}'_{g,r}(X,\beta)]^{vir}} 
\prod_{i=1}^r \psi_i^{k_i} \text{ev}_i^*(\gamma_{i}).$$
The associated partition function is defined 
by{\footnote{Our 
notation follows \cite{MNOP2,moop} and emphasizes the
role of the moduli space $\overline{M}'_{g,r}(X,\beta)$. 
The degree 0 collapsed contributions
will not appear anywhere in the paper.}} 
\begin{equation}
\label{abc}
\bZ'_{\mathsf{GW}}\Big(X;u\ \Big|\ \prod_{i=1}^r \tau_{k_i}(\gamma_{i})\Big)_\beta = 
\sum_{g\in{\mathbb Z}} \Big \langle \prod_{i=1}^r
\tau_{k_i}(\gamma_{i}) \Big \rangle'_{g,\beta} \ u^{2g-2}.
\end{equation}
Since the domain components must map nontrivially, an elementary
argument shows the genus $g$ in the  sum \eqref{abc} is bounded from below.

\subsection{Dimension constraints}
Descendents in Gromov-Witten and stable pairs theories are obtained via
tautological structures over the moduli spaces
$$\overline{M}'_{g,r}(X,\beta)\, , \ \ \ \ P_{n}(X,\beta)\times X$$
respectively.
The descendents $\tau_k(\gamma)$ in both cases mix the characteristic
classes of the tautological sheaves
$$\LL_i \rightarrow \overline{M}'_{g,r}(X,\beta)\, , \ \ \ \
\mathbb{F} \rightarrow P_{n}(X,\beta)\times X$$
 with the pull-back of
$\gamma\in H^*(X)$ via the evaluation/projective morphism.

In the absolute (nonequivariant) case, the Gromov-Witten and stable pairs
descendent series 
\begin{equation}\label{fhh6}
\bZ'_{\mathsf{GW}}\Big(X;u\ \Big|\ \prod_{i=1}^r \tau_{k_i}(\gamma_{i})\Big)_\beta\, , \ \ \ \
\bZ_{\mathsf{P}}\Big(X;q\ \Big|\ \prod_{i=1}^r \tau_{k_i}(\gamma_{i})\Big)_\beta
\end{equation}
both satisfy  dimension constraints.
For $\gamma_i\in H^{e_i}(X)$, the (real) dimension of the descendents 
Gromov-Witten and stable pairs theories are
$$\tau_{k_i}(\gamma_i)\in H^{{e_i}+2k_i}(\overline{M}'_{g,r}(X,\beta))\, , \ \ \ \
\tau_{k_i}(\gamma_i)\in H^{{e_i}+2k_i-2}(P_{n}(X,\beta))\, .$$
Since the virtual dimensions are
$$\text{dim}_{\mathbb{C}} \, [\overline{M}'_{g,r}(X,\beta)]^{vir} = \int_\beta c_1(T_X) + r\, , \ \ \ \
\text{dim}_{\mathbb{C}} \, [P_n(X,\beta)]^{vir} = \int_\beta c_1(T_X) $$
respectively, the dimension constraints
$$\sum_{i=1}^r \frac{e_i}{2}+k_i = \int_\beta c_1(T_X) + r \, , \ \ \ \
\sum_{i=1}^r \frac{e_i}{2}+k_i -1 = \int_\beta c_1(T_X) $$
exactly match.

After the matching of the dimension constraints, we can further reasonably ask if there is a relationship
between the Gromov-Witten and stable pairs descendent series \eqref{fhh6}.
The question has two immediately puzzling features:
\begin{enumerate}
\item[(i)] The series
involve different moduli spaces and universal structures.
\item[(ii)] The variables $u$ and $q$ of
the two series are different.
\end{enumerate}
Though the worry (i) is correct, both moduli spaces are
essentially based upon the geometry of curves in $X$, so there is hope
for a connection.
The {\em descendent correspondence}
proposes a precise relationship between the Gromov-Witten and stable pairs
descendent series, but only after a change of variables to
address (ii).

\subsection{Descendent notation}
%
Let $X$ be a nonsingular projective 3-fold.
Let $\widehat{\alpha}=(\widehat{\alpha}_1,\ldots,
\widehat{\alpha}_{\widehat{\ell}})$,
$$\widehat{\alpha}_1\geq \ldots\geq 
\widehat{\alpha}_{\widehat{\ell}} > 0\, ,$$
be a partition 
of size $|\widehat{\alpha}|$ 
and length $\widehat{\ell}$.
Let 
$$\iota_\Delta:\Delta\rightarrow X^{\widehat{\ell}}$$
 be the inclusion of the
small diagonal{\footnote{The small diagonal $\Delta$
is the set of points of $X^{\widehat{\ell}}$ for which
the coordinates $(x_1,\ldots, x_{\hat{\ell}})$ are all equal
$x_i=x_j$.}}
in the product $X^{\widehat{\ell}}$.
For $\gamma\in H^*(X)$, 
we write $$\gamma\cdot \Delta =\iota_{\Delta*}(\gamma) \in H^*(X^{\widehat{\ell}})\, .$$ 
Using the K\"unneth decomposition, we have
$$\gamma\cdot \Delta=
\sum_{{j_1, \ldots, j_{\hat{\ell}}}} c^\gamma_{j_1,\ldots, j_{\hat{\ell}}}\, 
\theta_{j_1} \otimes
\ldots\otimes \theta_{j_{\hat{\ell}}}\, ,$$
where $\{\theta_j\}$ is a $\mathbb{Q}$-basis of $H^*(X)$.
We define the descendent insertion $\tau_{\widehat{\alpha}}(\gamma)$ by
\begin{equation}\label{j77833}
\tau_{\widehat{\alpha}}(\gamma)= 
\sum_{j_1,\ldots,j_{\hat{\ell}}} c^\gamma_{j_1,\ldots, j_{\hat{\ell}}}\,
\tau_{\widehat{\alpha}_1-1}(\theta_{j_1})
\cdots\tau_{\widehat{\alpha}_{\hat{\ell}}-1}(\theta_{j_{\hat{\ell}}})\ .
\end{equation}
Three basic examples are:
\begin{enumerate}
\item[$\bullet$] If $\widehat{\alpha}=(\widehat{a}_1)$, then 
$$\tau_{(\, \widehat{a}_1\,)}(\gamma)= \tau_{\widehat{a}_1-1}(\gamma)\, .$$
The convention of shifting
the descendent by $1$  allows us to index descendent insertions by
standard partitions $\widehat{\alpha}$
 and follows the notation of \cite{PPDC}.

\item[$\bullet$] If $\widehat{\alpha}=(\widehat{a}_1,\widehat{a}_2)$ 
and $\gamma=1$ is the identity class, then
$$\tau_{(\, \widehat{a}_1,\, \widehat{a}_2\, )}(1)= \sum_{j_1,j_2} c^1_{j_1,j_2} 
\tau_{\widehat{a}_1-1}(\theta_{j_1})\, \tau_{\widehat{a}_2-1}(\theta_{j_2})\, ,$$
where $\Delta = \sum_{j_1,j_2} c^1_{j_1,j_2}\, \theta_{j_1} \otimes \theta_{j_2}$
is the standard K\"unneth decomposition of the diagonal in $X^2$.
\item[$\bullet$] If 
$\gamma$ is the class of a point, then 
\[
\tau_{\widehat{\alpha}}(\mathsf{p})=
\tau_{\widehat{\alpha}_1-1}(\mathsf{p})\cdots\tau_{\widehat{\alpha}_{\hat{\ell}}-1}(\mathsf{p}).
\]
\end{enumerate}
By the multilinearity of descendent insertions, 
formula  \eqref{j77833} does not depend upon the
basis choice $\{\theta_j\}$.


\subsection{Correspondence matrix} \label{corrmat}
 
A central result of \cite{PPDC} is
the construction of
a universal correspondence matrix $\widetilde{\mathsf{K}}$ 
indexed by partitions
$\alpha$ and $\widehat{\alpha}$ of positive size with{\footnote{Here, $i^2=-1$.}}
$$\widetilde{\mathsf{K}}_{\alpha,\widehat{\alpha}}\in 
\mathbb{Q}[i,c_1,c_2,c_3]((u))\, . $$
The elements of $\widetilde{\mathsf{K}}$
are constructed from the capped descendent vertex \cite{PPDC}
and satisfy two basic properties:

\begin{enumerate}
\item[(i)] The vanishing
$\widetilde{\mathsf{K}}_{\alpha,\widehat{\alpha}}=0$ holds {unless}
$|{\alpha}|\geq |\widehat{\alpha}|$.
\item[(ii)]
The $u$ coefficients of  $\widetilde{\mathsf{K}}_{\alpha,\widehat{\alpha}}\in
\mathbb{Q}[i,c_1,c_2,c_3]((u))$
are homogeneous{\footnote{The variable $c_i$ has degree $i$
for the homogeneity.}} in the variables $c_i$
of degree $$|\alpha|+\ell(\alpha) - |\widehat{\alpha}| 
- \ell(\widehat{\alpha})-3(\ell(\alpha)-1).$$ 
\end{enumerate}
Via the substitution
\begin{equation} \label{h3492}
c_i=c_i(T_X),
\end{equation}
the matrix elements  of $\widetilde{\mathsf{K}}$
act by cup product on the cohomology 
of $X$ with $\mathbb{Q}[i]((u))$-coefficients.

The matrix $\widetilde{\mathsf{K}}$ is 
used to define a correspondence
rule
\begin{equation}\label{pddff}
{\tau_{\alpha_1-1}(\gamma_1)\cdots
\tau_{\alpha_{\ell}-1}(\gamma_{\ell})}\ \  \mapsto\ \ 
\overline{\tau_{\alpha_1-1}(\gamma_1)\cdots
\tau_{\alpha_{\ell}-1}(\gamma_{\ell})}\ .
\end{equation}
The definition of the right side
of \eqref{pddff} requires a sum over all set
partitions $P$ of $\{ 1,\ldots, \ell \}$.
 For such a  set partition
$P$, each element $S\in P$
is a subset of $\{1,\ldots, \ell\}$.
Let $\alpha_S$ be the associated subpartition of
$\alpha$, and let
$$\gamma_S = \prod_{i\in S}\gamma_i.$$
In case all cohomology classes $\gamma_j$ are even,
we define the right side of the correspondence rule  \eqref{pddff} 
by
\begin{equation}\label{mqq23}
\overline{\tau_{\alpha_1-1}(\gamma_1)\cdots
\tau_{\alpha_{\ell}-1}(\gamma_{\ell})}
=
\sum_{P \text{ set partition of }\{1,\ldots,\ell\}}\ \prod_{S\in P}\ \sum_{\widehat{\alpha}}\tau_{\widehat{\alpha}}(\widetilde{\mathsf{K}}_{\alpha_S,\widehat{\alpha}}\cdot\gamma_S) \ .
\end{equation}
The second sum  in \eqref{mqq23} is over 
all partitions $\widehat{\alpha}$ of positive size. However, by the vanishing
of property (i),
$$\widetilde{\mathsf{K}}_{\alpha_S,\widehat{\alpha}}=0 \ \ \text{unless}
\ \ |{\alpha_S}|\geq |\widehat{\alpha}|\, , $$
the summation index  may be restricted to partitions $\widehat{\alpha}$ of positive size bounded
by $|\alpha_S|$.

Suppose  $|\alpha_S|=|\widehat{\alpha}|$
in the second sum in \eqref{mqq23}.
The homogeneity property (ii)
then places a strong constraint. The $u$ coefficients of  
\begin{equation*}
\widetilde{\mathsf{K}}_{\alpha_S,\widehat{\alpha}}\in
\mathbb{Q}[i,c_1,c_2,c_3]((u))
\end{equation*}
are homogeneous
of degree 
\begin{equation}\label{d2399}
3-2\ell(\alpha_S)  
- \ell(\widehat{\alpha}) \, .
\end{equation}
For the matrix element $\widetilde{\mathsf{K}}_{\alpha_S,\widehat{\alpha}}$ to be nonzero, the degree 
\eqref{d2399}
must be nonnegative. Since the lengths of $\alpha_S$ and 
$\widehat{\alpha}$ are at least 1,
nonnegativity of \eqref{d2399} is only possible if 
$$\ell(\alpha_S)= \ell(\widehat{\alpha})=1\, .$$
Then, we also have $\alpha_S=\widehat{\alpha}$ since the sizes match.

The above argument shows that the descendents on the right side of  \eqref{mqq23}
all correspond to partitions of size {\em less} than $|\alpha|$
except for the {\em leading term} obtained from the 
the maximal set partition
$$\{1\} \cup \{2\} \cup \ldots \cup \{\ell\} = \{1,2,\ldots, \ell\}$$
in $\ell$ parts.
The leading term of the descendent correspondence, calculated
in \cite{PPDC}, is a third basis property of $\widetilde{\mathsf{K}}$:

\begin{enumerate}
\item[(iii)]
$\ \ \overline{\tau_{\alpha_1-1}(\gamma_1)\cdots
\tau_{\alpha_{\ell}-1}(\gamma_{\ell})}
= (iu)^{\ell(\alpha)-|\alpha|}\, \tau_{\alpha_1-1}(\gamma_1)\cdots
\tau_{\alpha_{\ell}-1}(\gamma_{\ell}) +\ldots .$
\end{enumerate}


In case $\alpha=1^\ell$ has all parts equal to 1, then $\alpha_S$ also
has all parts equal to 1 for every $S\in P$. By property (ii), the $u$ coefficients of  
$\widetilde{\mathsf{K}}_{\alpha_S,\widehat{\alpha}}$
are homogeneous
of degree 
$$3-\ell(\alpha_S) - |\widehat{\alpha}| 
- \ell(\widehat{\alpha}),$$
and hence vanish unless
$$\alpha_S= \widehat{\alpha}=(1)\ .$$
Therefore, 
if
$\alpha$ has all parts equal to $1$,
the leading term is therefore the entire formula 
We obtain a fourth property of the matrix $\widetilde{\mathsf{K}}$:

\begin{enumerate}
\item[(iv)]
$\ \overline{\tau_{0}(\gamma_1)\cdots
\tau_{0}(\gamma_{\ell})}
=  \tau_{0}(\gamma_1)\cdots
\tau_{0}(\gamma_{\ell})\, .$
\end{enumerate}

In the presence of odd cohomology, a natural sign
must be included in formula \eqref{mqq23}. We may write 
set partitions $P$ of $\{1,\ldots, \ell\}$ indexing 
the sum on the right side of \eqref{mqq23} 
as
$$S_1\cup \ldots \cup S_{|P|} = \{1,\ldots, \ell\}.$$
The parts $S_i$ of $P$ are unordered, but we choose an ordering
for each $P$.
We then 
obtain a permutation of $\{1, \ldots, \ell\}$
by moving the elements to the ordered parts $S_i$ (and
respecting the original order in each group).
The permutation, in turn, determines a sign $\sigma(P)$
determined by the anti-commutation of the associated
odd classes. We then write
\begin{equation*}
\overline{\tau_{\alpha_1-1}(\gamma_1)\cdots
\tau_{\alpha_{\ell}-1}(\gamma_{\ell})}
=
\sum_{P \text{ set partition of }\{1,\ldots,\ell\}}\ (-1)^{\sigma(P)}
\prod_{S_i\in P}\ \sum_{\widehat{\alpha}}\tau_{\widehat{\alpha}}(\widetilde{\mathsf{K}}_{\alpha_{S_i},\widehat{\alpha}}
\cdot\gamma_{S_i}) \ .
\end{equation*}
The descendent 
$\overline{\tau_{\alpha_1-1}(\gamma_1)\cdots
\tau_{\alpha_{\ell}-1}(\gamma_{\ell})}$ is easily seen to have the
same commutation rules with respect to odd
cohomology as ${\tau_{\alpha_1-1}(\gamma_1)\cdots
\tau_{\alpha_{\ell}-1}(\gamma_{\ell})}$.

The geometric construction of $\widetilde{\mathsf{K}}$ in \cite{PPDC} 
expresses the coefficients explicitly in terms of the 1-legged capped
descendent vertex for stable pairs and stable maps.  These vertices
can be computed (as a rational function in the stable pairs case and
term by term in the genus parameter for stable maps). Hence, the coefficient
$$\widetilde{\mathsf{K}}_{\alpha,\widehat{\alpha}}\in \mathbb{Q}[i,c_1,c_2,c_3]((u))$$
can, in principle, be calculated term by term in $u$. The calculations
in practice are quite difficult, and 
 complete closed formulas are not known for all of the coefficients.

\subsection{Absolute case}
To state the descendent correspondence proposed in
\cite{PPDC} for all nonsingular projective 3-folds $X$, the basic degree
$$d_\beta = \int_{\beta} c_1(X) \ \in \mathbb{Z}$$
associated to the class $\beta\in H_2(X,\mathbb{Z})$ is required.

\begin{Conjecture}[P.-Pixton (2011)]
\label{ttt222} Let $X$ be a nonsingular projective 3-fold.
For $\gamma_i \in H^{*}(X)$, we have
\begin{multline*}
(-q)^{-d_\beta/2}\ZZ_{\mathsf{P}}\Big(X;q\ \Big|  
{\tau_{\alpha_1-1}(\gamma_1)\cdots
\tau_{\alpha_{\ell}-1}(\gamma_{\ell})}
\Big)_\beta \\ =
(-iu)^{d_\beta}\ZZ'_{\mathsf{GW}}\Big(X;u\ \Big|   
\ \overline{\tau_{\alpha_1-1}(\gamma_1)\cdots
\tau_{\alpha_{\ell}-1}(\gamma_{\ell})}\ 
\Big)_\beta 
\end{multline*}
under the variable change $-q=e^{iu}$.
\end{Conjecture}

Since the stable pairs side of the correspondence
 $$\ZZ_{\mathsf{P}}\Big(X;q\ \Big|  
{\tau_{\alpha_1-1}(\gamma_1)\cdots
\tau_{\alpha_{\ell}-1}(\gamma_{\ell})}
\Big)_\beta\, \in \mathbb{Q}((q))$$
 is  defined as a series in $q$, the change of variable
$-q=e^{iu}$ is {\em not} \`a priori well-defined.
However,
the stable pairs descendent series
is predicted by Conjecture \ref{111} to be a rational function in $q$.
The change of variable $-q=e^{iu}$ is well-defined for a rational function
in $q$ 
by substitution. The well-posedness of 
Conjecture \ref{ttt222} therefore depends upon Conjecture \ref{111}.

\subsection{Geometry of descendents}\label{ggogg}
Let $X$ be a nonsingular projective 3-fold, and let
 $D\subset X$ be a nonsingular divisor.
The Gromov-Witten descendent insertion $\tau_1(D)$ has a simple
geometric leading term.  
Let
$$[f:(C,p) \rightarrow X] \in \overline{M}_{g,1}(X,\beta)$$
be a stable map. Let
$$\text{ev}_1: \overline{M}_{g,1}(X,\beta) \rightarrow X$$
be the evaluation map at the marking.
The cycle
$$\text{ev}^{-1}_1(D) \subset \overline{M}_{g,1}(X,\beta)$$
corresponds to stable maps with $f(p)\in D$.
On the locus $\text{ev}^{-1}_1(D)$, there is a differential
\begin{equation}\label{sgg47}
df: T_{C,p} \rightarrow N_{X/D,f(p)}
\end{equation}
from the tangent space of $C$ at $p$ to the normal space
of $D\subset X$ at $f(p)\in D$. The differential
\eqref{sgg47} on $\text{ev}^{-1}_1(D)$ vanishes on the locus
where $f(C)$ is {\em tangent} to $D$ at $p$.
In other words,
$$\tau_{1}(D)+\tau_0(D^2) = \text{ev}_1^{-1}(D) \left( -c_1(T^*_{C,p})+
\text{ev}_1^*(N_{X/D})\right)$$
has the tangency cycle as a leading term. There are
correction terms from the loci where $p$ lies on a component of $C$ contracted by $f$ to
a point of $D$.

A parallel relationship can be pursued for $\tau_k(D)$ for for higher $k$
in terms of the locus of stable maps with higher tangency along $D$
at $f(p)$. A full correction calculus in 
case $X$ has dimension 1 (instead of 3) was found in \cite{OPGWH}.
The method has also been successfully applied to calculate the characteristic
numbers of curves in $\mathsf{P}^2$ for genus at most 2 in \cite{KGP}.{\footnote{In higher
genus, the correction calculus in $\mathsf{P}^2$ was too complicated
to easily control.}} 

By the Gromov-Witten/Pairs correspondence of Conjecture \ref{ttt222},
the stable pairs descendent $\tau_k(D)$ has leading term on the
Gromov-Witten side
$$\overline{\tau_k(D)} = (iu)^{-k} \tau_{k}(D) + \ldots \, .$$
Hence, the descendents $\tau_k(D)$ on the stable pairs side
should be viewed as essentially connected to the tangency loci
associated to the divisor $D\subset X$.

\subsection{Equivariant case}\label{eee999}
If $X$ is a nonsingular quasi-projective toric 3-fold,
all terms of the descendent correspondence
have $\T$-equivariant interpretations.
We take the equivariant K\"unneth decomposition in \eqref{j77833},
and the equivariant Chern classes $c_i(T_X)$ with respect to the
canonical $\T$-action on $T_X$ in \eqref{h3492}.
The toric case is proven in \cite{PPDC}.

\begin{Theorem}[P.-Pixton, 2011] Let \label{tt66}
$X$ be a nonsingular quasi-projective toric 3-fold.
For $\gamma_i \in H^{*}_{\mathbf{T}}(X)$,
we have
\begin{multline*}
(-q)^{-d_\beta/2}\ZZ_{\mathsf{P}}\Big(X;q\ \Big|  
{\tau_{\alpha_1-1}(\gamma_1)\cdots
\tau_{\alpha_{\ell}-1}(\gamma_{\ell})}
\Big)^\T_\beta \\ =
(-iu)^{d_\beta}\ZZ'_{\mathsf{GW}}\Big(X;u\ \Big|   \
\overline{\tau_{\alpha_1-1}(\gamma_1)\cdots
\tau_{\alpha_{\ell}-1}(\gamma_{\ell})}\
\Big)^\T_\beta 
\end{multline*}
under the variable change $-q=e^{iu}$. 
\end{Theorem}

Since the stable pairs side of the correspondence
 $$\ZZ_{\mathsf{P}}\Big(X;q\ \Big|  
{\tau_{\alpha_1-1}(\gamma_1)\cdots
\tau_{\alpha_{\ell}-1}(\gamma_{\ell})}
\Big)^\T_\beta\, \in \mathbb{Q}(s_1,s_2,s_3)((q))$$
 is a rational function in $q$ by 
Theorem \ref{pp12}, 
the change of variable
$-q=e^{iu}$ is well-defined by substitution.

When $X$ is a nonsingular projective
toric $3$-fold, Theorem  \ref{tt66}
implies Conjecture \ref{ttt222} for $X$ by taking the non-equivariant limit. However, Theorem \ref{tt66} is much 
stronger in the toric case than Conjecture \ref{ttt222} 
since the descendent insertions may exceed the virtual dimension in equivariant cohomology.

In case  $\alpha=(1)^\ell$ has all parts equal to 1, Theorem \ref{tt66}
specializes by property (iv) of Section \ref{corrmat} to the simpler statement
\begin{multline}\label{pp889}
(-q)^{-d_\beta/2}\ZZ_{\mathsf{P}}\Big(X;q\ \Big|\,  
{\tau_{0}(\gamma_1)\cdots
\tau_{0}(\gamma_{\ell})}
\Big)^\T_\beta \\ =
(-iu)^{d_\beta}\ZZ'_{\mathsf{GW}}\Big(X;u\ \Big|   \,
{\tau_{0}(\gamma_1)\cdots
\tau_{0}(\gamma_{\ell})}\,
\Big)^\T_\beta 
\end{multline}
which was first proven in the context of ideal sheaves in \cite{moop}.
Viewing both sides of \eqref{pp889} as series in $u$, we can complex conjugate
the coefficients. Imaginary numbers only occur in 
$$-q = e^{iu} \ \ \ \ \text{and} \ \ \ \ (-iu)^{d_\beta}\, .$$
After complex conjugation, we  find
\begin{multline*}
(-q)^{d_\beta/2}\ZZ_{\mathsf{P}}\Big(X;\frac{1}{q}\ \Big|  
\, {\tau_{0}(\gamma_1)\cdots
\tau_{0}(\gamma_{\ell})}\,
\Big)^\T_\beta \\ =
(iu)^{d_\beta}\ZZ'_{\mathsf{GW}}\Big(X;u\ \Big|   \,
{\tau_{0}(\gamma_1)\cdots
\tau_{0}(\gamma_{\ell})}\,
\Big)^\T_\beta 
\end{multline*}
and thus obtain the functional equation
\begin{equation*}
\ZZ_{\mathsf{P}}\Big(X;\frac{1}{q}\, \Big|\,  
{\tau_{0}(\gamma_1)\cdots
\tau_{0}(\gamma_{\ell})}\,
\Big)^\T_\beta  =
q^{-d_\beta}\ZZ_{\mathsf{P}}\Big(X;q\, \Big|   \,
{\tau_{0}(\gamma_1)\cdots
\tau_{0}(\gamma_{\ell})}\,
\Big)^\T_\beta 
\end{equation*}
as predicted by Conjecture \ref{444}.

\subsection{Relative case}
\subsubsection{Relative Gromov-Witten theory}
Let $X$ be a nonsingular projective 3-fold with a nonsingular divisor
$$D\subset X\, .$$ 
The relative theory of 
stable pairs was discussed in Section \ref{relth}.
A parallel relative Gromov-Witten theory of stable maps 
with specified tangency along
the divisor $D$ can also be defined.

In Gromov-Witten theory, relative conditions  are
represented by a partition $\mu$ of the integer
$
\int_\beta [D],
$
each part $\mu_i$ of which is marked
by a cohomology class $\delta_i\in H^*(D,\mathbb{Z})$,
\begin{equation}
\label{mm33}
\mu=( (\mu_1,\delta_1), \ldots, (\mu_\ell,\delta_\ell))\, .
\end{equation}
The numbers $\mu_i$ record
the multiplicities of intersection with $D$
while the cohomology labels $\delta_i$ record where the tangency occurs.
More precisely, let $\oM_{g,r}'(X/D,\beta)_\mu$ be the
moduli space of stable relative maps with tangency conditions
$\mu$ along $D$. To impose the full boundary condition,
we 
pull-back the classes $\delta_i$
via the evaluation maps
\begin{equation}\label{gtth341}
\oM_{g,r}'(X/D,\beta)_\mu \to D
\end{equation}
at the points of tangency.
Also, the tangency points are
considered to be unordered.{\footnote{The evaluation
maps are well-defined only after ordering the points.
We define the theory first with ordered tangency points. 
The unordered theory is then defined by dividing by the
automorphisms of the cohomology weighted partition $\mu$.}}

Relative Gromov-Witten theory was defined before the
study of stable pairs.
For the foundations, including the
definition of the moduli space of stable relative
maps and the construction of the virtual class
$$[\oM_{g,r}'(X/D,\beta)_\mu] \in H_*(\oM_{g,r}'(X/D,\beta)_\mu)\, ,$$
we refer the reader to \cite{LR,L}.

\subsubsection{Diagonal classes}\label{diagclas}
Definition \eqref{mqq23}  of the Gromov-Witten/Pairs correspondence
in the absolute case involves the diagonal 
$$\iota_\Delta:\Delta\rightarrow X^s$$
via \eqref{j77833}.
For the correspondence in the relative case, the diagonal
has a more subtle definition.

For the absolute geometry $X$, the product $X^s$ naturally
parameterizes $s$ ordered (possibly coincident) points on $X$.
For the relative geometry $X/D$, the parallel
object is the  moduli space $(X/D)^s$ of
$s$ ordered (possibly coincident) points 
$$(p_1,\ldots, p_s) \in X/D\, .$$
The points 
parameterized by $(X/D)^s$ 
are not allowed to lie on the relative divisor $D$.
When the points approach $D$, the target $X$ degenerates.
The resulting moduli space $(X/D)^s$ is a nonsingular variety.
Let
$$\Delta_{\mathsf{rel}} \subset (X/D)^s$$
be the small diagonal where all the points $p_i$ 
are coincident.
As a variety, $\Delta_{\mathsf{rel}}$ is isomorphic to $X$.

The space $(X/D)^s$
is a special case of well-known constructions in
relative geometry. 
For example, $(X/D)^2$ consists of
6 strata:

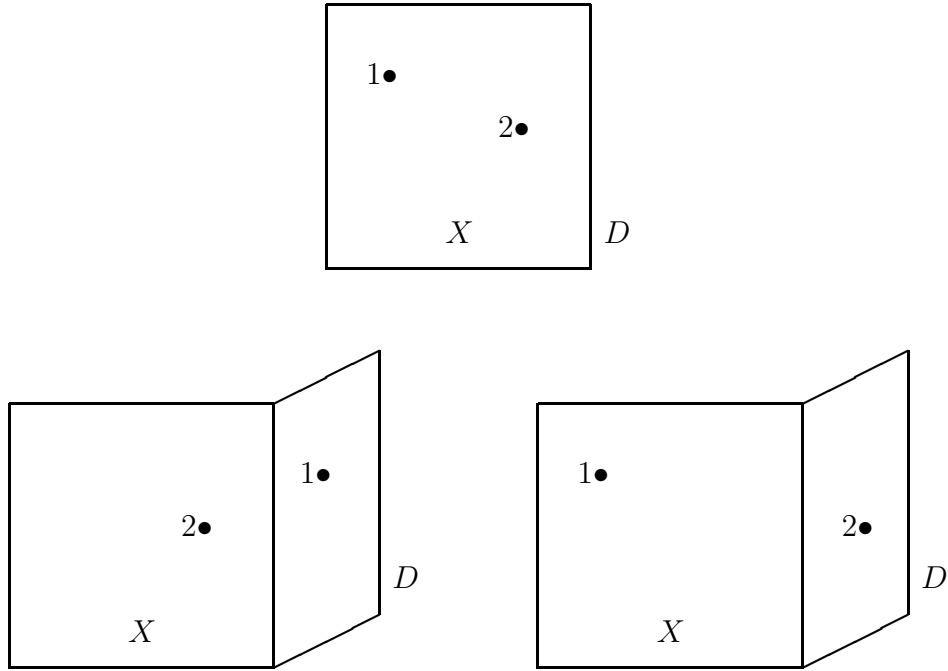
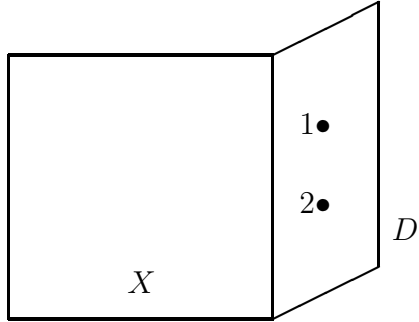
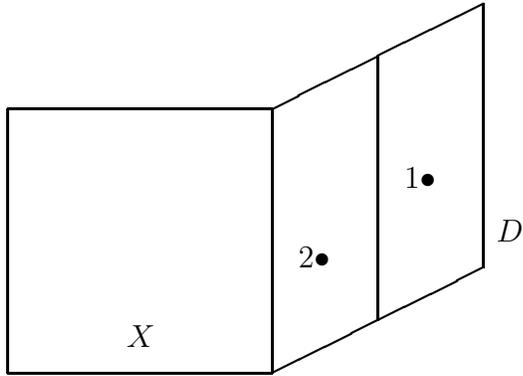
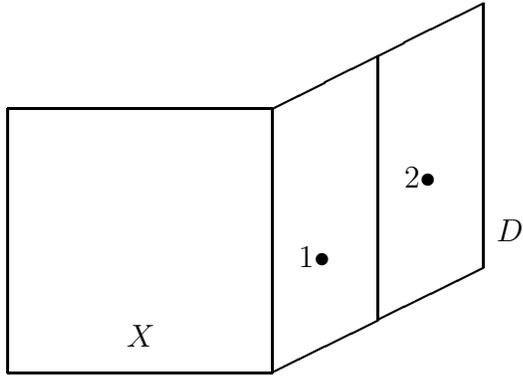
\begin{picture}(150,150)(-120,-5)
\thicklines
\put(10,10){\line(1,0){100}}
\put(10,110){\line(1,0){100}}
\put(10,10){\line(0,1){100}}
\put(110,10){\line(0,1){100}}
\put(25,80){$1\bullet$}
\put(75,60){$2\bullet$}
\put(55,20){$X$}
\put(115,20){$D$}
\end{picture}

\begin{picture}(150,150)(0,-5)
\thicklines
\put(10,10){\line(1,0){100}}
\put(10,110){\line(1,0){100}}
\put(10,10){\line(0,1){100}}
\put(110,10){\line(0,1){100}}

\put(110,110){\line(2,1){40}}
\put(110,10){\line(2,1){40}}
\put(150,30){\line(0,1){100}}
\put(155,40){$D$}

\put(120,80){$1\bullet$}
\put(75,60){$2\bullet$}
\put(55,20){$X$}

\put(210,10){\line(1,0){100}}
\put(210,110){\line(1,0){100}}
\put(210,10){\line(0,1){100}}
\put(310,10){\line(0,1){100}}

\put(310,110){\line(2,1){40}}
\put(310,10){\line(2,1){40}}
\put(350,30){\line(0,1){100}}
\put(355,40){$D$}

\put(225,80){$1\bullet$}
\put(325,60){$2\bullet$}
\put(255,20){$X$}

\end{picture}

\begin{picture}(150,150)(-100,-5)
\thicklines
\put(10,10){\line(1,0){100}}
\put(10,110){\line(1,0){100}}
\put(10,10){\line(0,1){100}}
\put(110,10){\line(0,1){100}}

\put(110,110){\line(2,1){40}}
\put(110,10){\line(2,1){40}}
\put(150,30){\line(0,1){100}}
\put(155,40){$D$}

\put(120,80){$1\bullet$}
\put(120,50){$2\bullet$}
\put(55,20){$X$}

\end{picture}

\begin{picture}(150,150)(-80,-5)
\thicklines
\put(10,10){\line(1,0){100}}
\put(10,110){\line(1,0){100}}
\put(10,10){\line(0,1){100}}
\put(110,10){\line(0,1){100}}

\put(110,110){\line(2,1){40}}
\put(110,10){\line(2,1){40}}
\put(150,30){\line(0,1){100}}

\put(150,30){\line(2,1){40}}
\put(150,130){\line(2,1){40}}
\put(190,50){\line(0,1){100}}
\put(195,60){$D$}

\put(160,80){$1\bullet$}
\put(120,50){$2\bullet$}
\put(55,20){$X$}

\thicklines

\end{picture}

\begin{picture}(150,150)(-80,-5)
\thicklines
\put(10,10){\line(1,0){100}}
\put(10,110){\line(1,0){100}}
\put(10,10){\line(0,1){100}}
\put(110,10){\line(0,1){100}}

\put(110,110){\line(2,1){40}}
\put(110,10){\line(2,1){40}}
\put(150,30){\line(0,1){100}}

\put(150,30){\line(2,1){40}}
\put(150,130){\line(2,1){40}}
\put(190,50){\line(0,1){100}}
\put(195,60){$D$}

\put(160,80){$2\bullet$}
\put(120,50){$1\bullet$}
\put(55,20){$X$}

\thicklines

\end{picture}

\noindent As a variety, $(X/D)^2$ is the blow-up of $X^2$ along $D^2$.
And, $\Delta_{\mathsf{rel}} \subset (X/D)^2$ is the strict transform
of the standard diagonal.

Select a subset $S$ of cardinality $s$
from the $r$ markings of the moduli space
of maps.
Just as $\oM_{g,r}'(X,\beta)$
admits a canonical evaluation to $X^s$ via
the selected markings, 
the moduli space $\oM_{g,r}'(X/D,\beta)_\mu$
admits a canonical evaluation
 $$\text{ev}_S: \oM_{g,r}'(X/D,\beta)_\mu \rightarrow (X/D)^s , $$
well-defined by the definition of a relative stable
map (the markings never map to the relative divisor).
The class 
$$\text{ev}_S^*(\Delta_{\mathsf{rel}}) \in H^*(\oM_{g,r}'(X/D,\beta)_\mu)$$
plays a crucial role in the relative descendent
correspondence.

By forgetting the relative structure, we obtain a projection
$$\pi:(X/D)^s \rightarrow X^s\ .$$
The product contains the standard diagonal $\Delta\subset X^s$. However,
$$\pi^*(\Delta) \neq \Delta_{\mathsf{rel}}\ .$$
The former has more components in the relative boundary if
$D\neq \emptyset$.

\subsubsection{Relative descendent correspondence} \label{pwwf}
Let $\widehat{\alpha}$ be a partition
of length $\widehat{\ell}$.
Let $\Delta_{\mathsf{rel}}$ be the cohomology class of the
small diagonal in  $(X/D)^{\widehat{\ell}}$.
For a cohomology class $\gamma$ of $X$, let
$$\gamma\cdot \Delta_{\mathsf{rel}} \in H^*\big((X/D)^{\widehat{\ell}}\, \big)\, ,$$
where $\Delta_{\mathsf{ref}}$ is the small diagonal of Section \ref{diagclas}.
Define the relative descendent insertion $\tau_{\widehat{\alpha}}(\gamma)$ by
\begin{equation}\label{j9994}
\tau_{\widehat{\alpha}}(\gamma)= 
\psi_1^{\widehat{\alpha}_1-1} \cdots \psi_{\hat{\ell}}^{\widehat{\alpha}_{\hat{\ell}}-1} \cdot
\text{ev}^*_{1,\ldots,\hat{\ell}} ( \gamma\cdot \Delta_{\mathsf{rel}}) \ .
\end{equation}
In case, $D=\emptyset$, definition \eqref{j9994}
specializes to \eqref{j77833}.

Let $\Omega_X[D]$ denote the locally free sheaf of 
differentials with logarithmic poles along $D$.
Let $$T_{X}[-D] = \Omega_{X}[D]^{\ \vee}$$
denote the dual sheaf of tangent fields with logarithmic
zeros.

For the relative geometry $X/D$, the coefficients of
the correspondence matrix 
$\widetilde{\mathsf{K}}$ act on the cohomology of $X$ via the
substitution
$$c_i= c_i(T_{X}[-D])$$
instead of the substitution $c_i=c_i(T_X)$ used in the absolute case.
Then, we define 
\begin{equation} \label{gtte4}
\overline{\tau_{\alpha_1-1}(\gamma_1)\cdots
\tau_{\alpha_{\ell}-1}(\gamma_{\ell})}
=
\sum_{P \text{ set partition of }\{1,\ldots,l\}}\ \prod_{S\in P}\ \sum_{\widehat{\alpha}}\tau_{\widehat{\alpha}}(\widetilde{\mathsf{K}}_{\alpha_S,\widehat{\alpha}}\cdot\gamma_S) \ 
\end{equation}
as before via \eqref{j9994} instead of 
\eqref{j77833}.
Definition \eqref{gtte4} is for  even classes $\gamma_i$.
In the presence of odd $\gamma_i$, a sign has to be included
exactly as in the absolute case.

\begin{Conjecture}
\label{ttt444} 
For $\gamma_i \in H^{*}(X)$,
we have 
\begin{multline*}
(-q)^{-d_\beta/2}\ZZ_{\mathsf{P}}\Big(X/D;q\ \Big|  
{\tau_{\alpha_1-1}(\gamma_1)\cdots
\tau_{\alpha_{\ell}-1}(\gamma_{\ell})} \ \Big| \ \mu
\Big)_\beta \\ =
(-iu)^{d_\beta+\ell(\mu)-|\mu|}\ZZ'_{\mathsf{GW}}\Big(X/D;u\ \Big|   
\ \overline{\tau_{a_1-1}(\gamma_1)\cdots
\tau_{\alpha_{\ell}-1}(\gamma_{\ell})}
\ \Big| \ \mu\Big)_\beta 
\end{multline*}
under the variable change $-q=e^{iu}$.
\end{Conjecture}

The change of variables is well-defined by the rationality of 
Conjecture \ref{222}.
A case in which Conjecture \ref{ttt444} is proven is when
$X$ is a nonsingular projective toric 3-fold and $D\subset X$ is
a toric divisor. The rationality of the stable pairs series
is given by Theorem \ref{ppp12}. The following result
can be obtained by the methods of \cite{PPQ}.

\begin{Theorem} 
For $X/D$ a nonsingular projective relative toric 3-fold,
 the descendent partition function
For $\gamma_i \in H^{*}(X)$,
we have 
\begin{multline*}
(-q)^{-d_\beta/2}\ZZ_{\mathsf{P}}\Big(X/D;q\ \Big|  
{\tau_{\alpha_1-1}(\gamma_1)\cdots
\tau_{\alpha_{\ell}-1}(\gamma_{\ell})} \ \Big| \ \mu
\Big)_\beta \\ =
(-iu)^{d_\beta+\ell(\mu)-|\mu|}\ZZ'_{\mathsf{GW}}\Big(X/D;u\ \Big|   
\ \overline{\tau_{a_1-1}(\gamma_1)\cdots
\tau_{\alpha_{\ell}-1}(\gamma_{\ell})}
\ \Big| \ \mu\Big)_\beta 
\end{multline*}
under the variable change $-q=e^{iu}$.
\end{Theorem}

Conjecture \ref{ttt444} can be lifted in a canonical
way to the equivariant relative case (as in the 
the rationality of Conjecture \ref{333}). Some 
equivariant relative results are proven in \cite{PPQ}.

\subsection{Complete intersections}
Let $X$ be a Fano or Calabi-Yau complete intersection of ample divisors
in a product of projective spaces,
$$ X \subset \PP^{n_1} \times \cdots \times \PP^{n_m}\ .$$
A central result of \cite{PPQ} is the proof of the descendent
correspondence for even classes.

\begin{Theorem} [P.-Pixton, 2012]
\label{qqq111} 
Let $X$ be a nonsingular Fano or Calabi-Yau complete intersection 3-fold in
a product of projective spaces. For even classes
$\gamma_i \in H^{2*}(X)$, we have
\begin{multline*}
(-q)^{-d_\beta/2}\ZZ_{\mathsf{P}}\Big(X;q\ \Big|  
{\tau_{\alpha_1-1}(\gamma_1)\cdots
\tau_{\alpha_{\ell}-1}(\gamma_{\ell})}
\Big)_\beta \\ =
(-iu)^{d_\beta}\ZZ'_{\mathsf{GW}}\Big(X;u\ \Big|   
\ \overline{\tau_{\alpha_1-1}(\gamma_1)\cdots
\tau_{\alpha_{\ell}-1}(\gamma_{\ell})}\ 
\Big)_\beta 
\end{multline*}
under the variable change $-q=e^{iu}$.
\end{Theorem}

Theorem \ref{qqq111} relies on the rationality of the stable
pairs series of 
Theorem \ref{qqq111f}. For
$\gamma_i \in H^{2*}(X)$ even classes of {\em positive} degree, 
we obtain from Theorem \ref{qqq111} (under the same complete intersection hypothesis for $X$) the following result where
only the leading term of the correspondence contributes:
\begin{multline*}
(-q)^{-d_\beta/2}\ \bZ_{\mathsf{P}}\left(X;q \
\Bigg| \ \prod_{i=1}^r {\tau}_0(\gamma_{i})
 \prod_{j=1}^s {\tau}_{k_j}(\mathsf{p}) \right)_{\beta}=\\
(-iu)^{d_\beta} (iu)^{-\sum k_j}\ 
 \bZ'_{\mathsf{GW}}\left(X;u \ \Bigg| \ \prod_{i=1}^r \tau_0(\gamma_{i}) 
\prod_{j=1}^s {\tau}_{k_j}(\mathsf{p})  \right)_{\beta} \ 
\end{multline*}
under the variable change $-q=e^{iu}$.

If we specialize Theorem \ref{qqq111} further 
to the case where there are no descendent insertions, we obtain 
$$
\ZZ_{\mathsf{P}}\Big(X;q\Big)_\beta  =
\ZZ'_{\mathsf{GW}}\Big(X;u\Big)_\beta 
$$
under the variable change $-q=e^{iu}$ for Calabi-Yau complete
intersections in a product of projective spaces.
In particular, the Gromov-Witten/Pairs correspondence
hold for the famous quintic Calabi-Yau 3-fold
$$X_5 \subset \proj^4\, .$$

\subsection{$K3$ fibrations}

Let $Y$ be a nonsingular projective toric 3-fold for which 
the anticanonical class $K_Y^*$ is base point free and the generic 
anticanonical
divisor is a nonsingular projective $K3$ surface $S$.
Let 
\begin{equation}\label{xhxhxh}
X \subset Y \times \mathsf{P}^1
\end{equation}
be a nonsingular hypersurface in the class $K^*_Y \otimes K^*_{\mathsf{P}^1}$.
Using the degeneration 
$$ X \leadsto Y\,  \cup\,  S\times \mathsf{P}^1\,  \cup\, Y$$
obtained by factoring a divisor of $K^*_Y \otimes K^*_{\mathsf{P}^1}$,
the results of \cite{PPQ} yield the Gromov-Witten/Pairs correspondence
for the Calabi-Yau 3-fold $X$.{\footnote{The strategy here is simpler
than presented in Appendix B of \cite{rp13} for a particular toric 4-fold
$Y$.}}

The hypersurface $X$ defined by \eqref{xhxhxh} is a $K3$-fibered
Calabi-Yau 3-fold. A very natural question to ask is whether
the Gromov-Witten/Pairs correspondence can be proven for all
$K3$-fibered 3-folds. While the general case is open, 
results for the correspondence in fiber classes
can be found in \cite{rp13}.{\footnote{Parallel questions can be pursued
for other surfaces. For results for surfaces of general type (involving
the stable pairs theory of descendents), see
\cite{KoolT}.}}

\section{Virasoro constraints} \label{333r}

\subsection{Overview}
Descendent partition functions in Gromov-Witten theory are
conjectured to satisfy Virasoro constraints \cite{EHX} for every target variety $X$.
Via the Gromov-Witten/Pairs descendent correspondence, 
we expect parallel constraints for the descendent theory of stable
pairs. An ideal path to finding the constraints for stable pairs
would be to start with the explicit Virasoro constraints in
Gromov-Witten theory and then apply the correspondence. However, our
knowledge of the correspondence matrix is not yet sufficient for such
an application.

Another method is to look experimentally for relations which are
of the expected shape. In a search conducted almost 10 years ago with A. 
Oblomkov and
A. Okounkov, we found a set of such relations for the theory of
ideal sheaves \cite{oop} for every nonsingular projective 3-fold $X$.
As an example, the equations for $\mathsf{P}^3$ are presented here
for stable pairs.{\footnote{Since
\cite{oop} is written for ideal sheaves, a DT/PT correspondence for descendents
is needed to move the relations to the theory of stable pairs. Such a correspondence
is also studied in \cite{oop}. I am very grateful to A. Oblomkov for 
his help with the formulas here.}}

\subsection{First equations}\label{firsteq}
Let $X$ be a nonsingular projective 3-fold.
The descendent insertions 
$$\tau_0(1)\, , \ \ \ \tau_0(D)\, \ \text{for $D\in H^2(X)$}, \ \ \ \tau_1(1)$$
all satisfy simple equations (parallel to the string, divisor, and dilation
equations in Gromov-Witten theory):
\begin{enumerate}
\item[(i)]
$\ZZ_{P}\Big(X;q\ \Big|\,  \tau_0(1)\cdot \prod_{i=1}^r \tau_{k_i}(\gamma_i)
\Big)_\beta = 0$,
\item[(ii)]
$\ZZ_{P}\Big(X;q\ \Big|\,  \tau_0(D)\cdot \prod_{i=1}^r \tau_{k_i}(\gamma_i)
\Big)_\beta = 
\left( \int_\beta D\right) \, \ZZ_{P}\Big(X;q\ \Big|  \prod_{i=1}^r \tau_{k_i}(\gamma_i)
\Big)_\beta \, $,
\item[(iii)]  $\ZZ_{P}\Big(X;q\ \Big|\,  \tau_1(1)\cdot \prod_{i=1}^r \tau_{k_i}(\gamma_i)
\Big)_\beta = 
\left(q\frac{d}{dq} - \frac{d_\beta}{2} \right) \, \ZZ_{P}\Big(X;q\ \Big|  \prod_{i=1}^r \tau_{k_i}(\gamma_i)
\Big)_\beta $ .
\end{enumerate}
All three are obtained directly from the definition of the descendent action given
 in Section \ref{actact}.
To prove (iii), the Hirzebruch-Riemann-Roch equation
$$\text{ch}_3(F)= n - \frac{d_\beta}{2}$$
is used for a stable pair
$$[F,s] \in P_n(X,\beta)\, , \ \ \ d_\beta= \int_\beta c_1(X)\, .$$

The compatibility of (i) and (ii) with the functional equation of Conjecture \ref{444} is trivial.
While not as obvious,  
the differential operator
$$q\frac{d}{dq} - \frac{d_\beta}{2}$$
is also beautifully consistent with  Conjecture \ref{444}.
We can easily prove using (iii) that Conjecture \ref{444} holds for
$$\ZZ_{P}\Big(X;q\ \Big|\,  \tau_1(1)\cdot \prod_{i=1}^r \tau_{k_i}(\gamma_i)
\Big)_\beta$$
if and only if 
Conjecture \ref{444} holds for 
$$\ZZ_{P}\Big(X;q\ \Big|\, \prod_{i=1}^r \tau_{k_i}(\gamma_i)
\Big)_\beta\, .$$
For example, equation (iii) yields
$$\ZZ_{\mathsf{P}}\big(\PP^3;q\ |\,   \tau_1(1) \tau_5(\mathsf{1}) \big
)_{\mathsf{L}} =  \frac{q+4q^2+17q^3-62q^4+17q^5+4q^6+q^7}{9(1+q)^4}\, $$
when applied to \eqref{s555}.

\subsection{Operators and constraints}
\label{sec:virasoro-constraints}
A basis of the cohomology $H^*(\mathsf{P}^3)$ is given by
$$\mathsf{1}\, ,\ \mathsf{H}\, ,\ \mathsf{L}=\mathsf{H}^2\, ,\ \mathsf{p}=\mathsf{H}^3$$
where $\mathsf{H}$ is the hyperplane class.
The divisor and dilaton equations here are
\begin{eqnarray*}
\ZZ_{P}\Big(\mathsf{P}^3;q\ \Big|\,  \tau_0(\mathsf{H})\cdot \mathsf{D})
\Big)_{d\mathsf{L}} &=& 
d \ZZ_{P}\Big(\mathsf{P}^3;q\ \Big|\,  \mathsf{D}\Big)_{d\mathsf{L}} \, , \\
\ZZ_{P}\Big(\mathsf{P}^3;q\ \Big|\,  \tau_1(1)\cdot \mathsf{D}\Big)_{d\mathsf{L}} &=& 
\left(q\frac{d}{dq} - 2{d} \right) \, \ZZ_{P}\Big(\mathsf{P}^3;q\ \Big|\, 
\mathsf{D}\Big)_{d\mathsf{L}} \, ,
\end{eqnarray*}
where $\mathsf{D}= \prod_{i=1}^r \tau_{k_i}(\gamma_i)$ is an arbitrary
descendent insertion.

Before presenting the formulas, we introduce two conventions
which simplify the notation.
The first concerns descendents with negative subscripts. We
define the descendent action in two negative cases:
\begin{equation}\label{jww9}
\tau_{-2}(\mathsf{H}^j)=-\delta_{j,3}\,,\quad \tau_{-1}(\gamma)=0\, .
\end{equation}
In particular, these all vanish except for $\tau_{-2}(\mathsf{p})= -1$.
Convention \eqref{jww9} is consistent with Definition \ref{dact} via the
replacement
$$\text{ch}_{2+i}(\mathbb{F}) \mapsto \text{ch}_{2+i}(\mathbb{I}[1]^\bullet)\, \, ,$$
where $\mathbb{I}^\bullet$ is the universal stable pair on
$X \times P_n(X,\beta)$.

For the Virasoro constraints, the formulas are more naturally stated
in terms of the Chern character subscripts (instead of including
the shift by 2 in Definition \ref{dact}). As a second convention, we
define the insertions $\mathsf{ch}_i(\gamma)$
by
\begin{equation}\label{kqq2}
\mathsf{ch}_i(\gamma)=\tau_{i-2}(\gamma) 
\end{equation}
for all $i\geq 0$. In particular, $\mathsf{ch}_0(\mathsf{p})$ acts
as $-1$ and $\mathsf{ch}_1(\mathsf{H}^j)$ acts as 0.

 Let $\mathbb{D}^+$ be a $\mathbb{Q}$-polynomial ring with generators
$$\Big\{\, \mathsf{ch}_i(\mathsf{H}^j)\, \Big| \  i\ge 0\, ,\ \ j=0,1,2,3\, \Big\}\, .$$
Via equation \eqref{kqq2}, we view $\mathbb{D}^+$ as an
extension
$$\mathbb{D} \subset \mathbb{D}^+\, $$
of the algebra of descendents defined in Section \ref{actactt}.
We define  
$$\mathsf{ch}_a\mathsf{ch}_b(\mathsf{H}^j) \in \mathbb{D}^+$$
in terms of the generators
by 
 $$\mathsf{ch}_a\mathsf{ch}_b(\mathsf{H}^j) = \sum_{r,s} \mathsf{ch}_a(\gamma^L_r)
\mathsf{ch}_b(\gamma^R_s)$$
where the sum is indexed by the K\"unneth decomposition 
$$\mathsf{H}^j\cdot \Delta = \sum_{r,s} \gamma^L_r \otimes \gamma^R_s \in H^*(\mathsf{P}^3
\times \mathsf{P}^3)\, $$
and $\Delta \subset \mathsf{P}^3
\times \mathsf{P}^3$ is the diagonal.
Both $\mathsf{ch}_i(\mathsf{H}^j)$ 
and $\mathsf{ch}_a\mathsf{ch}_b(\mathsf{H}^j)$
define operators on
$\mathbb{D}^+$ by multiplication.

To write the Virasoro relations, we will define
 derivations 
$$\mathrm{R}_k: \mathbb{D}^+ \rightarrow \mathbb{D}^+$$
for $k \geq -1$ by the following
 action on the generators of $\mathbb{D}^+$,
$$\mathrm{R}_k\left(\mathsf{ch}_i(\mathsf{H}^j)\right)=
\left(\, \prod_{n=0}^{k} (i+j-3+n)\, \right)\, 
\mathsf{ch}_{k+i}(\mathsf{H}^j)\, .$$
In case $k=-1$, the product on the right is empty and
$$\mathrm{R}_{-1}\left(\mathsf{ch}_i(\mathsf{H}^j)\right)=
\mathsf{ch}_{i-1}(\mathsf{H}^j)\, .$$

\begin{definition} Let $\mathcal{L}_k:\mathbb{D}^+\rightarrow
\mathbb{D}^+$
for $k\geq -1$ be the operator
\begin{eqnarray*}
 \mathcal{L}_k&=&-2\sum_{a+b=k+2}(-1)^{d^L d^R}(a+d^L-3)!(b+d^R-3)!\,
\mathsf{ch}_a\mathsf{ch}_b(\mathsf{H})\\ & & +
\sum_{a+b=k}a!b!\,\mathsf{ch}_a\mathsf{ch}_b(\mathsf{p})\\& &
+\, \mathrm{R}_k+(k+1)!\, \mathrm{R}_{-1}\mathsf{ch}_{k+1}(\mathsf{p})\, .
\end{eqnarray*}
\end{definition}

The first term in the formula for $\mathcal{L}_k$ requires
explanation. By definition,
\begin{equation}\label{ppqq22}
\mathsf{ch}_a\mathsf{ch}_b(\mathsf{H}) = \mathsf{ch}_a(\mathsf{p})
\mathsf{ch}_b(\mathsf{H})
+\mathsf{ch}_a(\mathsf{L})
\mathsf{ch}_b(\mathsf{L})
+ \mathsf{ch}_a(\mathsf{H})
\mathsf{ch}_b(\mathsf{p})
\end{equation}
via the three terms of the K\"unneth decomposition of $\mathsf{H}\cdot\Delta$.
The notation
$$(-1)^{d^L d^R}(a+d^L-3)!(b+d^R-3)!\,
\mathsf{ch}_a\mathsf{ch}_b(\mathsf{H})\, $$
is shorthand for the sum
\begin{eqnarray*}
& & (-1)^{3\cdot 1}(a+3-3)!(b+1-3)!\, \mathsf{ch}_a(\mathsf{p})
\mathsf{ch}_b(\mathsf{H})\\
&+&
(-1)^{2\cdot 2}(a+2-3)!(b+2-3)!\, 
\mathsf{ch}_a(\mathsf{L})
\mathsf{ch}_b(\mathsf{L})\\
&+& (-1)^{1\cdot 3}
(a+1-3)!(b+3-3)!\, 
\mathsf{ch}_a(\mathsf{H})
\mathsf{ch}_b(\mathsf{p})\, .
\end{eqnarray*}
The three summands of \eqref{ppqq22} are each weighted by the factor
$$(-1)^{d^L d^R}(a+d^L-3)!(b+d^R-3)!$$ where $d^L$ is the (complex) degree
of $\gamma^L$ and $d^R$ is the (complex) degree of $\gamma^R$ with
respect to the K\"unneth summand $\gamma^L\otimes \gamma^R$. 


In the second term of the formula, $a!b!\, 
\mathsf{ch}_a
\mathsf{ch}_b(\mathsf{p})$ can be expanded as
$$a!b!\, \mathsf{ch}_a\mathsf{ch}_b(\mathsf{p}) =a!b!\, \mathsf{ch}_a(\mathsf{p})
\mathsf{ch}_b(\mathsf{p})\, .$$
The summations over $a$ and $b$ in the first two terms in 
the formula for $\mathcal{L}_k$ require
$a\geq0$ and $b\geq 0$.  All factorials with negative
arguments vanish. 

For example, the formula for the first operator $\mathcal{L}_{-1}$ is
\begin{eqnarray*}
\mathcal{L}_{-1}&=& \mathsf{R}_{-1} + 0! \, \mathsf{R}_{-1} \mathsf{ch}_0(\mathsf{p}) \, .
\end{eqnarray*}
For $\mathcal{L}_{0}$, we have
\begin{eqnarray*}
\mathcal{L}_{0}&=& -2\cdot(-1)^{3\cdot 1}(0+3-3)!(2+1-3)!\, \mathsf{ch}_0(\mathsf{p})\mathsf{ch}_2(\mathsf{H})\\
       & & -2\cdot(-1)^{2\cdot 2}(1+2-3)!(1+2-3)!\, \mathsf{ch}_1(\mathsf{L})\mathsf{ch}_1(\mathsf{L})\\
 & & -2\cdot(-1)^{1\cdot 3}(2+1-3)!(0+3-3)!\, \mathsf{ch}_2(\mathsf{H})\mathsf{ch}_0(\mathsf{p})\\
& & +\mathsf{ch}_0(\mathsf{p})\mathsf{ch}_0(\mathsf{p})\\
&& +\mathsf{R}_{0} + \mathsf{R}_{-1} \mathsf{ch}_1(\mathsf{p})\, . 
\end{eqnarray*}
After simplification, we obtain
$$\mathcal{L}_0= 4 \mathsf{ch}_0(\mathsf{p})\mathsf{ch}_2(\mathsf{H})-2 \mathsf{ch}_1(\mathsf{L})\mathsf{ch}_1(\mathsf{L})
+\mathsf{ch}_0(\mathsf{p})\mathsf{ch}_0(\mathsf{p})
 +\mathsf{R}_{0} + \mathsf{R}_{-1} \mathsf{ch}_1(\mathsf{p})\, . 
$$
The operators $\mathcal{L}_k$ on $\mathbb{D}^+$ are conjectured to be
the analogs for stable pairs of the Virasoro constraints for the Gromov-Witten
theory of $\mathsf{P}^3$.

\begin{Conjecture}[Oblomkov-Okounkov-P.] \label{fpp55} 
We have
$$\mathsf{Z}_{\mathsf{P}}(\mathsf{P}^3;q \ | \, \mathcal{L}_k \, \mathsf{D})_{d\mathsf{L}}=0$$
for all  $k\geq -1$, for all
$\mathsf{D}\in\mathbb{D}^+$, and 
 for all curve classes $d\mathsf{L}$.
\end{Conjecture}

For example, for $k=-1$, Conjecture \ref{fpp55} states
$$\mathsf{Z}_{\mathsf{P}}(\mathsf{P}^3;q \ | \,\mathcal{L}_{-1} \mathsf{D})_{d\mathsf{L}}=0\, .$$
By the above calculation of $\mathcal{L}_{-1}$, 
\begin{eqnarray*}
\mathsf{Z}_{\mathsf{P}}(\mathsf{P}^3;q \ | \,\mathcal{L}_{-1} \mathsf{D})_{d\mathsf{L}}  & =&
\mathsf{Z}_{\mathsf{P}}\Big(\mathsf{P}^3;q \ \Big| \,(\mathsf{R}_{-1} + 0! \, \mathsf{R}_{-1} \mathsf{ch}_0(\mathsf{p}))\, \mathsf{D} \Big)_{d\mathsf{L}}\\
& = & \mathsf{Z}_{\mathsf{P}}\Big(
\mathsf{P}^3;q \ \Big| \,
(\mathsf{R}_{-1} - \mathsf{R}_{-1})\, \mathsf{D} \Big)_{d\mathsf{L}}\\
& = & 0\, ,
\end{eqnarray*}
where we have also used the descendent action $\mathsf{ch}_0(\mathsf{p})=-1$.
The claim
$$\mathsf{Z}_{\mathsf{P}}(\mathsf{P}^3;q \ | \,
\mathcal{L}_{0} \mathsf{D})_{d\mathsf{L}}=0\, .$$
is easily reduced to the divisor equation (ii) of Section \ref{firsteq}
and is also true. 

The first nontrivial assertion of 
Conjecture \ref{fpp55} occurs for $k=1$,
$$\mathsf{Z}_{\mathsf{P}}(\mathsf{P}^3;q \ | \,\mathcal{L}_{1} \mathsf{D})_{d\mathsf{L}}  =
\mathsf{Z}_{\mathsf{P}}\Big(\mathsf{P}^3;q \ \Big| \,\big(-4\mathsf{ch}_3(\mathsf{H}) +
\mathsf{R}_1 + 2\mathsf{ch}_2(\mathsf{p}) \mathsf{R}_{-1}\big)\, \mathsf{D} \Big)_{d\mathsf{L}}=0\, , $$
which is at the moment unproven.
For example, let $\mathsf{D}= \mathsf{ch}_3(\mathsf{p})$ and $d=1$.
We obtain a prediction for descendent series for $\mathsf{P}^3$,
$$
-4\mathsf{Z}_{\mathsf{P}}(\mathsf{ch}_3(\mathsf{H}) \mathsf{ch}_3(\mathsf{p})
)_{\mathsf{L}}
+12\mathsf{Z}_{\mathsf{P}}(\mathsf{ch}_4(\mathsf{p}))_{\mathsf{L}}
+2\mathsf{Z}_{\mathsf{P}}(\mathsf{ch}_2(\mathsf{p}) \mathsf{ch}_2(\mathsf{p})
)_{\mathsf{L}} =0\, ,$$
which can be checked  using the evaluations
\begin{align*}
\mathsf{Z}_{\mathsf{P}}(\mathsf{ch}_3(\mathsf{H}) \mathsf{ch}_3(\mathsf{p})
)_{\mathsf{L}}
& =&   
\mathsf{Z}_{\mathsf{P}}(\tau_1(\mathsf{H}) \tau_1(\mathsf{p})
)_{\mathsf{L}} &=&
\frac{3}{4}q - \frac{3}{2}q^2 +\frac{3}{4}q^3\, ,\\
\mathsf{Z}_{\mathsf{P}}(\mathsf{ch}_4(\mathsf{p}))_{\mathsf{L}} & =& 
\mathsf{Z}_{\mathsf{P}}(\tau_2(\mathsf{p}))_{\mathsf{L}} &=&
 \frac{1}{12}q
-\frac{5}{6}q^2 +\frac{1}{12}q^3\, , \\
\mathsf{Z}_{\mathsf{P}}(\mathsf{ch}_2(\mathsf{p}) \mathsf{ch}_2(\mathsf{p})
)_{\mathsf{L}} 
& = &
\mathsf{Z}_{\mathsf{P}}(\tau_0(\mathsf{p}) \tau_0(\mathsf{p})
)_{\mathsf{L}}  &=&
 q+ 2q^2+q^3 \, .
\end{align*}

\subsection{The bracket}
To find the Virasoro bracket, we introduce the operators
\begin{eqnarray*}
 L_k&=&
-2\sum_{a+b=k+2}(-1)^{d^L d^R}(a+d^L-3)!(b+d^R-3)!
\mathsf{ch}_a\mathsf{ch}_b(\mathsf{H})\\
& &+\sum_{a+b=k}a!b!\mathsf{ch}_a\mathsf{ch}_b(\mathsf{p})\\
& &+ \mathsf{R}_k\, .
\end{eqnarray*}
We then obtain the Virasoro relations and the bracket with $\mathsf{ch}_k(\mathsf{p})$,
$$[L_k,L_m]=(m-k)L_{k+m},\quad [L_n,k!\mathsf{ch}_k(\mathsf{p})]=k\cdot 
(k+n)!\mathsf{ch}_{n+k}(\mathsf{p}).$$
The operators $\mathcal{L}_k$ are expressed in terms of  $L_k$ by:
$$\mathcal{L}_k=L_k+(k+1)!L_{-1}\mathsf{ch}_{k+1}(\mathsf{p}).$$


\section{Virtual class in algebraic cobordism} \label{4448}
\subsection{Overview}
Let $X$ be nonsingular projective 3-fold.
From the work of J. Shen \cite{Shen}, the virtual fundamental class
of the moduli space of stable pairs
$$[P_n(X,\beta)]^{vir} \in A_{d_\beta}(P_n(X,\beta))$$
admits a canonical lift to the theory of algebraic cobordism{\footnote{We
not do review the foundations of the theory of algebraic
cobordism here. 
The reader can find discussions in \cite{LMo, LPa}.
As for cohomology, we always take $\mathbb{Q}$-coefficients.
Shen constructs a canonical lift to algebraic cobordism 
$[M]^{vir} \in \Omega_*(M)$
of the
virtual class in Chow $[M]^{vir} \in A_*(M)$
obtained from a 2-term perfect obstruction theory on
a quasi-projective scheme $M$.
}} 
\begin{equation}\label{f99f}
[P_n(X,\beta)]^{vir} \in \Omega_{d_\beta}(P_n(X,\beta))\, 
\end{equation}
where $d_\beta=\int_\beta c_1(X)$.
Shen's construction depends only upon the 2-term perfect
obstruction theory of $P_n(X,\beta)$ and is
closely related to
earlier work of Ciocan-Fontantine and Kapranov \cite{CFK} and
Lowrey-Sch\"urg \cite{LS}. 

The lift \eqref{f99f} leads to several natural questions.
The simplest is {\em how does the virtual class in algebraic 
cobordism vary with $n$?}
Let
$$\pi: P_n(X,\beta) \rightarrow \bullet$$
be the structure map to the point $\bullet$.
Then, for fixed $\beta$, we define
$$\mathsf{Z}^{\Omega}_{\mathsf{P}}(X;q)_\beta= \sum_{n\in \mathbb{Z}}q^n\,  \pi_*[P_n(X,\beta)]^{vir}
\ \in \Omega_{d_\beta}(\bullet)
\otimes_{\mathbb{Q}} \mathbb{Q}((q))\, .$$
Is there an analogue for $\mathsf{Z}^{\Omega}_{\mathsf{P}}(X;q)_\beta$
of the rationality and functional
equation for the descendent theory of the
standard virtual class?

\subsection{Chern numbers}
While the full data of the cobordism class \eqref{f99f} is 
difficult to analyze, the push-forward
$$\pi_*[P_n(X,\beta)]^{vir} \in \Omega_{d_\beta}(\bullet)$$
is characterized by the virtual Chern numbers of
$P_n(X,\beta)$.

Since $P_n(X,\beta)$ has a 2-term perfect obstruction theory, there
is a virtual tangent complex
$\mathsf{T}^{vir} \in D^b(P_n(X,\beta))$
with Chern classes
$$c_i(\mathsf{T}^{vir})\in H^{2i}( P_n(X,\beta))\, .$$
For every partition of the virtual dimension $d_\beta$,
$$\sigma=(s_1,\ldots, s_\ell)\, , \ \ \ \ \ d_\beta=\sum_{i=1}^\ell s_i\, ,$$
we define an associated Chern number
$$c^\sigma_{n,\beta} = \int_{[P_n(X,\beta)]^{vir}} \prod_{i=1}^\ell c_{s_i}(\mathsf{T}^{vir})\ \in 
\mathbb{Z}$$
by integration against the standard virtual class
$$[P_n(X,\beta)]^{vir} \in H_{2d_\beta}(P_n(X,\beta))\, .$$
The complete collection of Chern numbers
$$\big\{\, c^\sigma_n \, \big| \, \sigma \in \text{Partitions}(d_\beta)\, \big\}$$
uniquely determines the algebraic cobordism class
$$\pi_*[P_n(X,\beta)]^{vir} \in \Omega_{d_\beta}(\bullet)\, .$$

\subsection{Rationality and the functional equation}
 The rationality of the partition function $\mathsf{Z}^{\Omega}_{\mathsf{P}}(X;q)_\beta$
is equivalent to the rationality of {\em all} the functions
$$\mathsf{Z}^\sigma_{\mathsf{P}}(X;q)_\beta = \sum_{n\in \mathbb{Z}} c^\sigma_{n,\beta} q^n$$
for  $\sigma \in \text{Partitions}(d_\beta)$.

\begin{Theorem}[Shen 2014]\label{ll22}
The Chern class $c_i(\mathsf{T}^{vir})\in H^{2i}( P_n(X,\beta))$ can be written as
a $\mathbb{Q}$-linear combination of products of descendent classes
$$\left\{ \, \prod_{i=1}^r \tau_{k_i}(\gamma_i) \ \Big| \ \sum_{i=1}^rk_i \equiv 0 \, \text{\em mod 2}\, , \, \gamma_i\in H^*(X)\, \right\}$$
by a formula which is independent of $n$ and $\beta$.
\end{Theorem} 

Shen's proof is geometric and constructive. Following the notation of Section \ref{actact},
let $$\pi_P: X \times P_n(X,\beta) \rightarrow P_n(X,\beta)$$
be the projection and let $\mathbb{I}^\bullet\in D^b( X \times P_n(X,\beta) )$
be the universal stable pair.
The class of the virtual tangent complex in 
$K^0(P_n(X,\beta))$ is 
\begin{eqnarray*}
[-\mathsf{T}^{vir}] &=& [R\pi_{P*} R\mathcal{H}om(\mathbb{I}^\bullet,\mathbb{I}^\bullet)_0] \\
& = & [R\pi_{P*} (\mathbb{I}^\bullet\otimes^L (\mathbb{I}^\bullet))^\vee]
- [R\pi_{P*} \mathcal{O}_{X \times P_n(X,\beta)}]\, .
\end{eqnarray*}
The Chern character of $-\mathsf{T}^{vir}$ is then computed by the Grothendieck-Riemman-Roch
formula,
\begin{equation}\label{nana4}
\text{ch}[-\mathsf{T}^{vir}] = \pi_{P*}\Big(\text{ch}(\mathbb{I}^\bullet) \cdot \text{ch}((\mathbb{I}^\bullet)^\vee)
\cdot \text{Td}(X)\Big)
- \pi_{P*}\Big( \text{Td(X)}\Big)\, .
\end{equation}
The second term of\eqref{nana4} is just $\int_X \text{Td}_3(X)$
times the identity $1\in H^0(P_n(X,\beta))$.

More interesting is the first term of \eqref{nana4} which can be written
as 
\begin{equation}\label{jj33}
\epsilon_{*}\Big(\text{ch}(\mathbb{I}^\bullet) \cdot \text{ch}((\widetilde{\mathbb{I}}^\bullet)^\vee)
\cdot \Delta \cdot \text{Td}(X)\Big)
\end{equation}
where $\epsilon$ is the projection
$$\epsilon: X \times X \times P_{n}(X,\beta)\rightarrow P_n(X,\beta)\, ,$$
$\mathbb{I}^\bullet$ and $\widetilde{\mathbb{I}}^\bullet$ are the universal
stable pairs pulled-back via the first and second projections
$$ X\times P_n(X,\beta)\, \leftarrow\, X \times X \times P_{n}(X,\beta)
\, \rightarrow\,  X\times P_n(X,\beta)$$
respectively, and
$\Delta$ is the pull-back of the diagonal in $X\times X$.
Using the K\"unneth decomposition of $\Delta$, Shen easily writes \eqref{jj33}
as a quadratic expression in the descendent classes --- see \cite[Section 3.1]{Shen}.
The answer is a universal formula independent of $n$ and $\beta$.

Though not explicitly remarked (nor needed) in \cite{Shen}, Shen's universal
formula for $\text{ch}[-\mathsf{T}^{vir}]$
is  a $\mathbb{Q}$-linear combination of classes
$$\left\{ \, \tau_{k_1}(\gamma_1)\tau_{k_2}(\gamma_2)\, \Big| \,  k_1+k_2\equiv 0 \, \text{ mod 2}\, , \ \gamma_1,\gamma_2 \in H^*(X) \, \right\}
$$
since each quadratic term appears in \eqref{jj33} in a form proportional to 
$$ ((-1)^{k_1}+(-1)^{k_2})\cdot \tau_{k_1}(\gamma_1)\tau_{k_2}(\gamma_2)$$
because of the universal stable pair 
$\text{ch}(\mathbb{I}^\bullet)$ appears together with 
the dual $\text{ch}((\widetilde{\mathbb{I}}^\bullet)^\vee)$.

There are two immediate consequences of Theorem \ref{ll22}. 
If the rationality of descendent series of Conjecture \ref{111}
holds for $X$, then
$${\text{\em $\mathsf{Z}^\Omega_{\mathsf{P}}(X;q)_\beta$ is the Laurent expansion of 
a rational function in $\Omega_{d_\beta}(\bullet)\otimes_{\mathbb{Q}}\mathbb{Q}(q)$}}\, .$$
In particular, Shen's results yield the rationality of the partition
functions in algebraic cobordism in case $X$ is a nonsingular projective
toric variety (where rationality of the descendent series is proven).

The second consequence concerns the functional equation.
The descendents which arise in Theorem \ref{ll22} have {\em even} 
subscript sum. Hence, if the functional equation of Conjecture \ref{444}
holds for $X$, then
\begin{equation}\label{gtt99}
\mathsf{Z}^\Omega_{\mathsf{P}}\left(X;\frac{1}{q}\right)_\beta = q^{-d_\beta} \mathsf{Z}^\Omega_{\mathsf{P}}(X;q)_\beta\,.
\end{equation}
The functional equation \eqref{gtt99} should be regarded
as the correct generalization to all $X$ of the symmetry
\begin{equation*}
\mathsf{Z}_{\mathsf{P}}\left(Y;\frac{1}{q}\right)_\beta = 
 \mathsf{Z}_{\mathsf{P}}(Y;q)_\beta\,
\end{equation*}
of stable pairs invariants for {\em Calabi-Yau} 3-folds $Y$.

\subsection{An example}
A geometric basis of $\Omega_*(\bullet)$ is given by the classes of products of projective spaces.
As an example, we write the series
$$\mathsf{Z}^\Omega_{\mathsf{P}}(\mathsf{P}^3;q)_{\mathsf{L}} \in \Omega_4(\bullet) \otimes_{\mathbb{Q}} 
\mathbb{Q}(q)$$
in terms of products of projective spaces:
\begin{eqnarray*}
\mathsf{Z}^\Omega_{\mathsf{P}}(\mathsf{P}^3;q)_{\mathsf{L}}&=&  
 \hspace{9pt} [\mathbb{P}^4]\cdot f_4(q)\\
&& + [\mathbb{P}^3 \times \mathbb{P}^1]\cdot f_{31}(q) \\
&& +[\mathbb{P}^2 \times \mathbb{P}^2]\cdot f_{22}(q)\\
 && +[\mathbb{P}^2 \times \mathbb{P}^1 \times \mathbb{P}^1] \cdot  f_{211}(q)\\
&& + [\mathbb{P}^1 \times \mathbb{P}^1 \times \mathbb{P}^1 \times \mathbb{P}^1]\cdot f_{1111}(q) \, ,
\end{eqnarray*}
where the rational functions{\footnote{I am very grateful to J. Shen for providing
these formulas.}} are given by

{\footnotesize{
\begin{eqnarray*}
f_4(q)&=& -4q -40q^2-4q^3\,, \\
f_{31}(q)&=& \frac{q}{(1+q)^4}\left( \frac{21}{2} +139q+ \frac{823}{2}{q^2}+446q^3+  \frac{823}{2}{q^4}+139q^5+ \frac{21}{2}q^6 \right)\, ,\\
f_{22}(q)&=& 6q + 60q^2 +6q^3\, , \\
f_{211}(q)&=& \frac{q}{(1+q)^4}\left( -18 -264q -774{q^2} - 816q^3 -774{q^4} - 264q^5 -18q^6 \right)\, ,\\
 f_{1111}(q)&=& \frac{q}{(1+q)^6}\Big(\frac{13}{2} +115q +490 q^2  +889q^3 + 1215q^4\\
& & \hspace{130pt} +889q^5 +490q^6  + 115q^7 + \frac{13}{2}q^8\Big)\,.
   \end{eqnarray*} }}

\subsection{Further directions}
The study of the virtual class in algebraic cobordism of the moduli space of stable pairs
$P_n(X,\beta)$ is intimately connected with the study of descendents invariants. The basic reason is because the Chern classes of the virtual 
tangent complex
are {\em tautological classes} of $P_n(X,\beta)$ in the sense of
 Section \ref{actactt}. 
If another approach to the virtual class in algebraic cobordism class could be found, perhaps
the implications could be reversed and results about descendent series could be proven.

\vspace{+4 pt}
\noindent Departement Mathematik \\
\noindent ETH Z\"urich \hfill  \\
\noindent rahul@math.ethz.ch

\end{document}